\documentclass[12pt]{article}

\usepackage{a4wide}
\usepackage{amssymb, amsmath}
\usepackage[all]{xypic}
\usepackage[mathscr]{eucal}

\newtheorem{theorem}[equation]{Theorem}
\newtheorem{proposition}[equation]{Proposition}
\newtheorem{lemma}[equation]{Lemma}

\newtheorem{corollary}[equation]{Corollary}

\renewcommand{\theequation}{\arabic{section}.\arabic{equation}}
\newcommand{\reset}{\setcounter{equation}{0}}

\newcommand{\Remark}{\addtocounter{equation}{1}\paragraph{Remark \theequation}}

\newcommand{\QQ}{\mathbb{Q}}
\newcommand{\NN}{\mathbb{N}}
\newcommand{\ZZ}{\mathbb{Z}}

\newcommand{\CC}{\mathbb{C}}
\renewcommand{\AA}{\mathbb{A}}
\newcommand{\PP}{\mathbb{P}}
\newcommand{\FF}{\mathbb{F}}

\newcommand{\M}{\mathbb{M}}

\newcommand{\fa}{\mathfrak{a}}

\newcommand{\fd}{\mathfrak{d}}
\newcommand{\ff}{\mathfrak{f}}
\newcommand{\fp}{\mathfrak{p}}
\newcommand{\fq}{\mathfrak{q}}
\newcommand{\fm}{\mathfrak{m}}
\newcommand{\fn}{\mathfrak{n}}

\newcommand{\fN}{\mathfrak{N}}

\newcommand{\fP}{\mathfrak{P}}

\newcommand{\MM}{\mathcal{M}}  

\newcommand{\OO}{\mathcal{O}}
\newcommand{\KK}{\mathcal{K}}
\newcommand{\TT}{\mathcal{T}}
\newcommand{\II}{\mathcal{I}}
\newcommand{\JJ}{\mathcal{J}}
\newcommand{\C}{\mathcal{C}}
\newcommand{\cQ}{\mathcal{Q}}

\newcommand{\Fq}{\FF_q}
\newcommand{\kinf}{k_{\infty}}
\newcommand{\Cinf}{\CC_{\infty}}
\newcommand{\Af}{\AA_f}        
\newcommand{\Ahat}{\hat{A}}

\newcommand{\End}{\mathrm{End}}
\newcommand{\Aut}{\mathrm{Aut}}
\newcommand{\Gal}{\mathrm{Gal}}

\newcommand{\sep}{\mathrm{sep}}
\newcommand{\ab}{\mathrm{ab}}
\newcommand{\an}{\mathrm{an}}
\newcommand{\Pic}{\mathrm{Pic}}
\newcommand{\Spec}{\mathrm{Spec}}
\newcommand{\diag}{\mathrm{diag}}

\newcommand{\Stab}{\mathrm{Stab}}
\newcommand{\Tr}{\mathrm{Tr}}

\newcommand{\GL}{\mathrm{GL}_2}

\newcommand{\PGL}{\mathrm{PGL}_2}
\newcommand{\PSL}{\mathrm{PSL}_2}
\newcommand{\PGLki}{\mathrm{PGL}_2(\kinf)}
\newcommand{\PSLki}{\mathrm{PSL}_2(\kinf)}

\newcommand{\HCM}{H_{\mathrm{CM}}}
\newcommand{\SC}{\mathcal{SC}}

\newcommand{\proof}{\noindent {\bf Proof. \;}}
\def\qed{{\hskip0pt\unskip\unskip\nobreak\hfil\penalty50
          \hskip1em\hbox{}\nobreak\hfil
           {$\square$}
          \parfillskip=0pt\finalhyphendemerits=0
          \par}\medskip}

\DeclareMathOperator*{\limproj}{\underleftarrow{\lim}}

\newbox\mybox
\def\arrover#1{\mathrel{
       \setbox\mybox=\hbox spread 1.4em
              {\hfil$\scriptstyle#1$\hfil}
       \vbox{\offinterlineskip\copy\mybox
             \hbox to\wd\mybox{\rightarrowfill}}}}

\begin{document}

\title{CM points on products of Drinfeld modular curves}
\author{Florian Breuer}
\maketitle

\begin{abstract}
Let $X$ be a product of Drinfeld modular curves over a general base ring $A$ of odd characteristic. We classify those  subvarieties of $X$ which contain a Zariski-dense subset of CM points. This is an analogue of the Andr\'e-Oort conjecture. As an application, we construct non-trivial
families of higher Heegner points on modular elliptic curves over
global function fields.
%
\end{abstract}

{\small

}

\section{Introduction}
\reset

\subsection{Basic notations}  
The following notations will be used throughout this paper. Let $\Fq$ denote
the finite field with $q$ elements, where $q$ is a power of the {\em odd} prime
$p$. Let $k$ be a global function field with field of constants $\Fq$. Fix a
place $\infty$ of $k$, and denote by $\kinf$ the completion of $k$ at $\infty$,
and by  $\Cinf = \hat{\bar{k}}_{\infty}$ the completion of an algebraic closure
of $\kinf$. Let $A = \{ x\in k \;|\; \text{$x$ is regular outside $\infty$}\}$, it
is a Dedekind domain of finite class number. We denote by $|\cdot|$ the absolute
value corresponding to $\infty$, and note that for $a\in A$, we have
$|a|=q^{\deg(a)}=|A/\langle a \rangle|$.

Unless stated otherwise, a {\em Drinfeld module} always means a rank $2$
Drinfeld $A$-module defined over $\Cinf$ of generic characteristic. See
\cite{Goss} and \cite{Hayes92} for an overview of Drinfeld modules.

We denote by $\Ahat = \limproj A/\fn$ the profinite completion of $A$, and by
$\Af = \Ahat \otimes_A k$ the ring of finite ad\`eles of $k$. Let $\KK$ be a
subgroup of finite index (i.e. an open subgroup) of $\GL(\Ahat)$, then we
denote by $\MM_A^2(\KK)$ the coarse moduli scheme parameterizing rank $2$
Drinfeld $A$-modules with level $\KK$-structure, it is an affine curve over 
$\Spec(A)$, which is not in general irreducible. A {\em Drinfeld modular curve}
over $\Cinf$ is an irreducible component of some  $\MM_A^2(\KK)_{\Cinf} =
\MM_A^2(\KK)\times_A\Cinf$. If $\KK=\GL(\Ahat)$ we will just write $\MM_A^2 =
\MM_A^2(\GL(\Ahat))$, which is the coarse moduli space of Drinfeld modules
without level structure.

Let $\Omega := \PP^1(\Cinf) - \PP^1(\kinf)$ denote the {\em Drinfeld upper
half-plane}, which has a rigid analytic structure. The group $\GL(\kinf)$ acts
on $\Omega$ via fractional linear transformations. Let $X$ be a Drinfeld modular
curve, then it is known that, as rigid analytic varieties, we have 
\begin{equation}\label{param}
X(\Cinf)^{\an} \cong \Gamma\backslash\Omega
\end{equation}
for some arithmetic subgroup $\Gamma \subset \GL(k)$. We denote by 
$[\omega]\in X(\Cinf)$ the point corresponding to $\omega\in\Omega$.

\subsection{Main results}

Let $X = X_1 \times \cdots \times X_n$ be a product of Drinfeld modular curves.
A point $x\in X(\Cinf)$ is called a {\em CM point} if it corresponds to an
$n$-tuple of Drinfeld modules with complex multiplication. These in turn
correspond via (\ref{param}) to points $\omega\in\Omega$ with $[k(\omega):k]=2$.

An irreducible algebraic subcurve $Y \subset X$ is called a {\em special subcurve} if there exists a
partition $\{1,\ldots,n\} = I_0 \coprod I_1$, $I_1 \neq\emptyset$, and
elements  $g_i \in \GL(k)$, for all $i\in I_1$, such that $Y = \{x\}\times Y'$,
where $x\in\prod_{i\in I_0}X_i(\Cinf)$ is a CM point and
$Y' \subset \prod_{i\in I_1}X_i$ is a curve such that $Y'(\Cinf)$ is the image 
of the map
\[ 
\Omega \longrightarrow \prod_{i\in I_1}X_i(\Cinf);
\quad \omega \longmapsto \big([g_i(\omega)]\big)_{i\in I_1}. 
\]
The curve $Y$ is called a {\em pure special subcurve} if $I_0 = \emptyset$,
i.e. if the projections $p_i : Y \rightarrow X_i$ are surjective for all
$i=1,\ldots,n$.

An irreducible algebraic subvariety $Y\subset X$ is called a {\em special subvariety} if there exists
a partition $\{1,\ldots,n\} = \coprod_{j=0}^g I_j$ such that  $Y = \{x\}\times
\prod_{j=1}^g Y_j$, where $x\in\prod_{i\in I_0}X_i(\Cinf)$ is a CM point, and
each $Y_j \subset \prod_{i\in I_j}X_i$ is a pure special subcurve, for
$j=1,\ldots,g$. We see that CM points are just special subvarieties of
dimension zero.

The aim of this paper is to prove the following result.

\begin{theorem}\label{MainResult}
Let $X = X_1\times\cdots\times X_n$ be a product of Drinfeld modular curves.
Then an irreducible algebraic subvariety $Y\subset X$ contains a Zariski-dense
subset of CM points if and only if $Y$ is a special subvariety.
\end{theorem}

It is easy to see that a special subvariety contains a Zariski-dense set of CM 
points, the hard part is to prove the converse. In the special case where $Y$
is a curve, we actually prove an effective result, namely that $Y$ is special
if and only if $Y$ contains a CM point of sufficiently large {\em CM height},
see Theorem \ref{MainCurves}.

Theorem \ref{MainResult} is an analogue of the Andr\'e-Oort Conjecture for
products of classical modular curves, see \cite{Andre} and \cite{Edixhoven}. In
fact, our proof is closely modeled on Edixhoven's approach \cite{Edixhoven}.

Theorem \ref{MainResult} was proved in \cite{BreuerPrep} for the special case
$A = \Fq[T]$. In this paper we show how to adapt the arguments of
\cite{BreuerPrep} to the case of general $A$ (but still of odd characteristic).
As an application, in Section \ref{HeegnerSection} we extend our previous results 
\cite{Breuer04} concerning higher
Heegner points on elliptic curves over rational functions fields to the case
of global function fields.

\section{Some preliminaries}
\reset

We begin by gathering some basic results which will be needed in the proof of
Theorem~\ref{MainResult}.

\subsection{Complex multiplication and CM heights}

Let $\phi$ be a Drinfeld module with complex multiplication by the ring $R$,
and let $K$ be the quotient field of $R$. Then $K/k$ is a quadratic {\em
imaginary} extension, which means that only one prime of $K$ lies above
$\infty$, which we again call $\infty$. Denote by $\OO_K$ the ring of integers
of $K$, i.e. the integral closure of $A$ in $K$. Then $A\subset R\subset\OO_K$,
and $R$ is a projective $A$-module of rank 2, hence by the invariant factor
theorem, $R=A+\ff\OO_K$, for some ideal $\ff\subset A$, which we will call the
conductor of $R$. Note that $\ff\OO_K$ is the largest $\OO_K$-ideal which is
also an $R$-ideal, which is the definition of conductor usually found in the
literature.

The ring $R$ is an order in $K$, and is not in general integrally closed.
However, one may still view $\phi$ as a Drinfeld $R$-module of rank $1$, after
Hayes \cite{Hayes79}, and in fact we have $\MM_R^1 = \Spec(\OO_{H_R})$, where
$\OO_{H_R}$ is the ring of integers of the class field $H_R/K$ corresponding to
the class group
$K^{\times}\backslash\AA_{f,K}^{\times}/\hat{R}^{\times}\cong\Pic(R)$.  The
field $H_R$ is also known as the {\em ring class field} of $R$, and $H_R/K$ is
unramified outside $\ff\OO_K$. Note that we only deal with the {\em finite}
ad\`eles $\AA_{f,K}$, so that $H_R/K$ splits completely at~$\infty$.  The
action of $\AA_{f,K}^{\times}$ on $\Spec(\OO_{H_R})$ coincides, 
via class field theory, with the action of $\Gal(H_R/K)$ on~$\MM_R^1$. Hence $\phi$ is
defined over $H_R$ and isogenies act like Galois. In particular, we have

\begin{proposition}\label{CMmain}
Let $\phi$ be a Drinfeld module with complex multiplication by $R$.
Let $\fn\subset A$ be a non-zero ideal such that every prime factor of $\fn$
splits in $R$, and let $\sigma_{\fn} = (\fn,H_R/k)\in\Gal(H_R/k)$ be the
corresponding Frobenius element. Then $\phi$ and $\phi^{\sigma_{\fn}}$ are 
linked by a cyclic isogeny of degree $\fn$.
\end{proposition}

\proof As every prime factor of $\fn$ splits in $R$, we may choose an ideal
$\fN$ of $R$ such that $R/\fN \cong A/\fn$. Now an isogeny of $\phi$ as a
rank $2$ Drinfeld $A$-module with kernel $A/\fn$ is also an isogeny of $\phi$ as
a rank $1$ Drinfeld $R$-module with kernel $R/\fN$. The result now follows from
the above discussion. \qed

Note that $H_R/k$ might not be abelian, but since all prime factors of $\fn$ split
in $K/k$ and $H_R/K$ is abelian, the conjugacy class $(\fn, H_R/K)$ contains only the one element
$\sigma_\fn = (\fN, H_R/K)$.

Denote by $|\ff|=|A/\ff|$, then we define the {\em CM height} of $\phi$ by
\begin{equation}
\HCM(\phi) := q^g|\ff|,
\end{equation}
where $g$ is the genus of $K$. Note that this definition differs from 
\cite[Def. 3.7]{BreuerPrep} by a power of $1/2$. If
$x=(x_1,\ldots,x_n)\in
\MM_A^2(\KK_1)(\Cinf)\times\cdots\times\MM_A^2(\KK_n)(\Cinf)$, then we define
\[
\HCM(x) := \max\big(\HCM(x_1),\ldots,\HCM(x_n)\big).
\]

The following result shows that $\HCM$ is a counting function on the set of CM
points in $\MM_A^2(\Cinf)$, which justifies the terminology.

\begin{proposition}\label{CMheight}
For every $B>0$, the set
\[
\{ \phi\in\MM_A^2(\Cinf) \;|\; \text{$\phi$ is CM and $\HCM(\phi)<B$}\}
\] 
is finite.
\end{proposition}

\proof For a given $g\geq 0$ there are only finitely many global function fields
$K$ with genus $g$ and field of constants contained in $\FF_{q^2}$. For each
such field $K$, there are only finitely many orders of the form $R=A+\ff\OO_K$
with bounded $\ff$. And for each such $R$, there are only $|\Pic(R)|$ Drinfeld
modules $\phi$ with $\End(\phi)\cong R$.\qed

We will need the following class-number estimate

\begin{proposition}\label{ClassNumber}
Let $\phi$ be a Drinfeld module with complex multiplication by $R$. Then for
every $\varepsilon>0$ there exists a computable constant $C_{\varepsilon}>0$
such that
\[
|\Pic(R)| > C_{\varepsilon}\HCM(\phi)^{1-\varepsilon}. 
\]
\end{proposition}

\proof Let $K$ be the quotient field of $R$, and $g$ its genus.
Firstly, we have \cite[Prop. 3.1]{BreuerPrep}
\[
|\Pic(\OO_K)| \geq h(K) \geq \frac{(q-1)(q^{2g}-2gq^g+1)}{2g(q^{g+1}-1)},
\]
where $h(K)$ denotes the class number of $K$.
Secondly, the exact sequence \cite[\S I.12]{Neukirch} 
\[
1\rightarrow\OO_K^{\times}/R^{\times}\longrightarrow
(\OO_K/\ff\OO_K)^{\times}/(R/\ff R)^{\times}\longrightarrow
\Pic(R)\longrightarrow\Pic(\OO_K)\rightarrow 1
\]
leads, as in the classical case, to the pleasing expression
\begin{equation}
|\Pic(R)| = \frac{|\Pic(\OO_K)|}{[\OO_K^{\times}:R^{\times}]}\cdot|\ff|
\prod_{\fp|\ff}\big(1-\chi(\fp)|\fp|^{-1}\big),
\end{equation}
where
\[
\chi(\fp) = \left\{\begin{array}{rl}
  1  & \text{if $\fp$ splits in $K/k$}\\
  -1 & \text{if $\fp$ is inert in $K/k$}\\
  0  & \text{if $\fp$ is ramified in $K/k$.}\end{array}\right.
\]
Here one uses the fact that $R/\ff R\cong A/\ff$. The estimate now follows
easily. \qed

\subsection{Drinfeld modular curves}

Let $\KK\subset\GL(\Ahat)$ be a subgroup of finite index, and recall that
$\MM_A^2(\KK)$ denotes the coarse moduli scheme parameterizing Drinfeld modules
with level-$\KK$ structure. Identifying Drinfeld modules over $\Cinf$ with
their associated rank $2$ lattices in $\Cinf$, and parameterizing the space of
such lattices ad\`elically, one arrives at the following analytic
parametrization (see e.g. \cite{GekelerDMC}).

\begin{eqnarray}
\MM_A^2(\KK)(\Cinf)^{\an} & \cong & 
  \GL(k)\backslash \Omega \times\GL(\Af)/\KK\nonumber\\
& \cong & \coprod_{s\in S}\Gamma_s\backslash\Omega,\label{thingy}
\end{eqnarray}
where $S\subset\GL(\Af)$ denotes a (finite) set of representatives for 
$\GL(k)\backslash\GL(\Af)/\KK$, and 
$\Gamma_s = s\KK s^{-1}\cap\GL(k)$. 

The determinant map induces an isomorphism
\[
\det : \GL(k)\backslash\GL(\Af)/\KK \stackrel{\sim}{\longrightarrow}
k^{\times}\backslash\Af^{\times}/\det(\KK)\quad\text{(\;$\cong \Pic(A)$ if
$\det(\KK)=\Ahat^{\times}$).}
\]

We next describe a scheme-theoretic version of the determinant map. Let $f\in
A$, and denote by $\KK(f)\subset\GL(\Ahat)$ the kernel of reduction mod $f$. 
The Weil pairing for Drinfeld modules \cite[Chap. 5]{vdHeidenThesis} gives us a
map of $A[1/f]$-schemes
\begin{equation}\label{WeilPairing}
\MM_A^2(\KK(f)) \stackrel{w_f}{\longrightarrow} \MM_A^1(\det(\KK(f))) = 
\Spec(\OO_{H_f}),
\end{equation}
where $\OO_{H_f}$ denotes the ring of integers of the class field $H_f$ of $k$
corresponding to the class group 
$k^{\times}\backslash\Af^{\times}/\det(\KK(f)) = 
k^{\times}\backslash\Af^{\times}/\{x\in \Ahat^{\times} \;|\; x \equiv 1 \bmod
f\}$. Under this map, the left $\Af^{\times}\GL(\Ahat)$-action on $\MM_A^2(\KK(f))$
corresponds to its determinant $\Af^{\times}$-action on $\Spec(\OO_{H_f})$,
which in turn corresponds to the $\Gal(k^{\ab}/k)$-action on both sides given
by the reciprocity map. It follows that $\MM_A^2(\KK(f))\times_A k$ is defined
over $k$, and all its $\bar{k}$-irreducible components are defined over $H_f$.

Suppose now that $\KK$ contains $\KK(f)$ for some $f\in A$, so $\KK$ is a
{\em congruence subgroup}. Then we may divide out by the action of $\KK$ and
$\det(\KK)$ in the left and right hand sides of (\ref{WeilPairing}),
respectively. As $\MM_A^2(\KK) = \KK\backslash\MM_A^2(\KK(f))$, we thus get a map of
$A[1/f]$-schemes
\begin{equation}\label{WeilPairing2}
\MM_A^2(\KK) \stackrel{w_{\KK}}{\longrightarrow} \MM_A^1(\det(\KK)) = 
\Spec(\OO_{H_{\KK}}),
\end{equation}
where this time $H_{\KK}$ is the class field corresponding to 
$k^{\times}\backslash\Af^{\times}/\det(\KK)$. As before, it follows that 
$\MM_A^2(\KK)\times_A k$ is defined over $k$, and its $\bar{k}$-irreducible
components are defined over $H_{\KK}$. 

We point out that, analytically, the maps $w_f$ and $w_{\KK}$ above send the
rank $2$ lattice $L$ to the rank $1$ lattice $\wedge^2 L$.

In the case $\KK=\GL(\Ahat)$, recall that we use the notation
$\MM_A^2:=\MM_A^2(\GL(\Ahat))$. Let $H$ denote the Hilbert class field of $k$,
corresponding to the class group
$k^{\times}\backslash\Af^{\times}/\Ahat^{\times}\cong\Pic(A)$. It is the
maximal unramified abelian extension of $k$ in which $\infty$ splits
completely. Then we see that the irreducible components of $\MM^2_{A,\Cinf}$
are defined over $H$. We also see that its set of irreducible components
corresponds to $\Pic(A)$ in  such a way that the component corresponding to
$[\fa]\in\Pic(A)$ parametrizes lattices isomorphic to $A\oplus\fa$ as
projective $A$-modules.

We may choose our representatives $S$ such that $1\in S$, in which case
$\Gamma_1 = \GL(A)$, and one of the irreducible components of $\MM^2_{A,\Cinf}$
thus corresponds to $\GL(A)\backslash\Omega$ via (\ref{thingy}). We call it the
{\em identity component}, and denote it by $\M$. It is the simplest modular
curve. If $A=\Fq[T]$ then in fact $\M \cong \AA^1$ is the affine line. In
general, we have seen that $\M$ is a smooth irreducible affine curve defined
over $H$.

\subsection{The curves $Y(\fn)$, $Y_0(\fn)$ and $Y_2(\fn)$}
\label{Y0section}

We denote by $Z$ the center of the algebraic group $\GL$. Let $\fn\subset A$
be a non-zero ideal, and consider the following open subgroups of
$\GL(\Ahat)$.
\begin{eqnarray*}
\KK(\fn)   & = & \ker\big(\GL(\Ahat) \rightarrow \GL(A/\fn)\big),\\
\KK_0(\fn) & = & \{\big(\begin{smallmatrix} a & b \\ c & d 
                 \end{smallmatrix}\big) \in \GL(\Ahat) \;|\; c\in\fn\Ahat\},\\
\KK_2(\fn) & = & \{\gamma\in\GL(\Ahat) \;|\; (\gamma \bmod \fn)\in 
                 Z(A/\fn)\}.
\end{eqnarray*}
These lead to coarse moduli schemes $\MM_A^2(\KK(\fn))$ (also denoted
$\MM_A^r(\fn)$ in the literature), $\MM_A^2(\KK_0(\fn))$ and
$\MM_A^2(\KK_2(\fn))$. Similarly to the previous section, we define the
Drinfeld modular curves $Y(\fn)$, $Y_0(\fn)$ and $Y_2(\fn)$ as the identity
components of the respective moduli schemes tensored with $\Cinf$.

We have isomorphisms of rigid analytic varieties
\[
Y(\fn)(\Cinf)^{\an}\cong\Gamma(\fn)\backslash\Omega,\quad 
Y_0(\fn)(\Cinf)^{\an}\cong\Gamma_0(\fn)\backslash\Omega,\quad 
Y_2(\fn)(\Cinf)^{\an}\cong\Gamma_2(\fn)\backslash\Omega,
\]
for the arithmetic groups
\begin{eqnarray*}
\Gamma(\fn)   & = & \{\gamma\in\GL(A) \;|\; \gamma \equiv 1 \bmod \fn\},\\
\Gamma_0(\fn) & = & \{\big(\begin{smallmatrix} a & b \\ c & d 
                      \end{smallmatrix}\big) \in \GL(A) \;|\; c\in\fn\},\\
\Gamma_2(\fn) & = & \{\gamma\in\GL(A) \;|\; (\gamma \bmod \fn)\in Z(A/\fn)\}.
\end{eqnarray*}
As in the classical case, $Y_0(\fn)$ parametrizes pairs of Drinfeld modules
linked by cyclic $\fn$-isogenies, but where the first Drinfeld module corresponds to a lattice isomorphic to $A^2$ as an $A$-module.

As $\KK(\fn)$, $\KK_0(\fn)$ and $\KK_2(\fn)$ are all congruence subgroups of
$\GL(\Ahat)$ (pick any $0\neq f\in\fn$), it follows from the Weil-pairing
(\ref{WeilPairing2}) that $Y_0(\fn)$ is defined over $H$, whereas the curves
$Y(\fn)$ and $Y_2(\fn)$ are defined over the class fields corresponding to
\[
k^{\times}\backslash\Af^{\times}/\{x\in\Ahat^{\times}\;|\;x\equiv 1 \bmod\fn\},
\quad\text{and}
\]
\[
k^{\times}\backslash\Af^{\times}/\{x\in\Ahat^{\times}\;|\;
\text{$x$ is a square mod $\fn$}\},
\]
respectively.

For a ring $R\supset\Fq$, and algebraic group $G$, we define 
\[
G^1(R) := \{g\in G(R) \;|\; \det(g) \in\Fq^{\times} \}.
\]
The {\em degree} of an ideal $\fn\subset A$ is 
$\deg(\fn):=\log_q|\fn|=\log_q|A/\fn|$. 

\begin{proposition}\label{covering}
We have covers $Y(\fn)\rightarrow Y_2(\fn)\rightarrow Y_0(\fn)\rightarrow\M$.
Moreover, $Y(\fn)$ and $Y_2(\fn)$ are Galois covers of $\M$ with Galois groups
\begin{eqnarray*}
\Gal(Y(\fn)/\M)   & \cong & \GL^1(A/\fn)/Z(\Fq) \\
\Gal(Y_2(\fn)/\M) & \cong & \GL^1(A/\fn)/Z^1(A/\fn)\cong\PGL^1(A/\fn). 
\end{eqnarray*}
Note that $\PGL^1(A/\fn)\cong\PSL(A/\fn)$ if every prime factor of $\fn$ has
even degree.
\end{proposition}

\proof Clearly, $\Gamma(\fn)\subset\Gamma_2(\fn)\subset\Gamma_0(\fn)$, whence
the coverings. Furthermore, $\Gamma(\fn)$ and $\Gamma_2(\fn)$ are normal
subgroups of $\GL(A)$, so the respective coverings are Galois with Galois groups
\[
\Gal(Y(\fn)/\M)\cong\GL(A)/\Gamma(\fn)\cong\GL^1(A/\fn)/Z(\Fq),
\quad\text{and}
\]
\[
\Gal(Y_2(\fn)/\M)\cong\GL(A)/\Gamma_2(\fn)\cong\GL^1(A/\fn)/Z^1(A/\fn)
\cong\PGL^1(A/\fn). 
\]
Lastly, if every prime factor of $\fn$ has even degree, then every element of
$\Fq$ is a square in $A/\fn$ and $\PGL^1(A/\fn)=\PSL(A/\fn)$.\qed

\subsection{Special subcurves of $\big(\MM^2_{A,\Cinf}\big)^n$}

In the Introduction we defined the notion of (pure) special subcurves of a
product of $n$ Drinfeld modular curves, in particular in $\M^n$.
 
For a subset $I\subset\{1,\ldots,m\}$ we denote by $p_I : \M^n \rightarrow \M^I$
the projection onto the coordinates listed in $I$. We will also write
$p_{\{i\}}=p_i$ and $p_{\{i,j\}}=p_{i,j}$.

The curve $Y_0(\fn)$ parametrizes Drinfeld modules linked by a cyclic isogeny of
degree $\fn$. If $(\phi,\phi')$ is such a pair, then $\phi\in\M(\Cinf)$; and
$\phi'\in\M(\Cinf)$ if and only if $\fn$ is a principal ideal. Thus, for
$\fn = \langle N \rangle$ principal, we obtain a map 
$Y_0(\fn) \rightarrow\M^2$ whose image is denoted by $Y'_0(N)$. Analytically, 
this map is given by
\begin{eqnarray}
\Gamma_0(\fn)\backslash\Omega & \longrightarrow & \nonumber
  \big(\GL(A)\backslash\Omega\big) \times \big(\GL(A)\backslash\Omega\big) \\
\;[\omega] & \longmapsto & \big([\omega],[N\omega]\big).\label{Y'}
\end{eqnarray}
In particular, we see that $Y'_0(N)$ is an irreducible pure special subcurve of
$\M^2$. Moreover, every irreducible pure special subcurve of $\M^2$ is of the
form $Y'_0(N)$ for some $N\in A$. Indeed, for any $g\in\GL(k)$, we denote by $N
= \deg(g)$ the determinant of $a\cdot g$, with $a\in A$ chosen in such a way
that the entries of $a\cdot g$ are in $A$ and have no common factor. This $N$
is defined up to the square of a unit, and is called the {\em degree} of $g$.
Now, let  $g_1,g_2 \in\GL(k)$ and $N:=\deg(g_2g_1^{-1})$. Then the special
curve corresponding to the map $\omega \mapsto
\big([g_1(\omega)],[g_2(\omega)]\big)$ is just $Y'_0(N)\subset\M^2$.

It follows that an irreducible curve $Y\subset\M^n$ is a pure special subcurve if
and only if $p_{i,j}(Y) = Y'_0(N_{ij})$ for some $N_{ij}\in A$ for all 
$i\neq j$.

Next, it will be convenient to study the space $\SC_n$ of all special subcurves
in $\big(\MM^2_{A,\Cinf}\big)^n$. Analytically, this will turn out to be the
following double coset space
\begin{equation}
\SC_n(\Cinf)^{\an} := 
 \GL(k)\backslash\GL(\Af)^n\times\Omega/\big(Z(\Af)\GL(\Ahat)\big)^n,
\end{equation}
where $\big(Z(\Af)\GL(\Ahat)\big)^n$ acts from the right on $\GL(\Af)^n$ in the usual
way, and trivially on $\Omega$, while $\GL(k)$ acts from the left on $\Omega$
in the usual way, and diagonally on $\GL(\Af)^n$.

Choose a set of representatives $T\subset\GL(\Af)^n$ for the double quotient 
\begin{equation}\label{cosetT}
\GL(k)\backslash\GL(\Af)^n/\big(Z(\Af)\GL(\Ahat)\big)^n.
\end{equation}
Note that $T$ is infinite, and we again choose $1\in T$. For each
$t=(t_1,\ldots,t_n)\in T$ let  $\Gamma_t = \cap_{i=1}^n
t_i\GL(\Ahat)t^{-1}_i\cap\GL(k)$. Then we have a canonical bijection 
\begin{eqnarray*}
\GL(k)\backslash\GL(\Af)^n\times\Omega/\big(Z(\Af)\GL(\Ahat)\big)^n
 & \stackrel{\sim}{\longrightarrow} & \coprod_{t\in T}\Gamma_t\backslash\Omega\\
\;[(h_1,\ldots,h_n),\omega]\hspace{2.4cm} & \longmapsto & [g^{-1}(\omega)]_{t},
\end{eqnarray*}
where $t\in T$ is such that $[t]=[(h_1,\ldots,h_n)]=[h]$ in (\ref{cosetT}), and $g\in\GL(k)$ is
such that $h=gta$, for some $a\in\big(Z(\Af)\GL(\Ahat)\big)^n$. We see that $\SC_n$ is
the disjoint union of an infinite family of Drinfeld modular curves,
parametrized by $T$.

Next, we describe a map from $\SC_n(\Cinf)$ into $\big(\MM^2_A(\Cinf)\big)^n$.
\begin{eqnarray}
\SC_n(\Cinf) & \stackrel{\theta}{\longrightarrow} & \big(\MM^2_A(\Cinf)\big)^n
\nonumber\\
\GL(k)\backslash\GL(\Af)^n\times\Omega/\big(Z(\Af)\GL(\Ahat)\big)^n & \longrightarrow &
\GL(k)^n\backslash\GL(\Af)^n\times\Omega^n/\big(Z(\Af)\GL(\Ahat)\big)^n \nonumber\\
& & = \big(\GL(k)\backslash\GL(\Af)\times\Omega/\GL(\Ahat)\big)^n\nonumber\\
\;[(h_1,\ldots,h_n),\omega] & \longmapsto & \big([h_1,\omega],\ldots,
[h_n,\omega]\big);
\nonumber\\
\; & \; & \;\nonumber\\ 
\coprod_{t\in T}\Gamma_t\backslash\Omega & \longrightarrow &
\coprod_{s\in S^n}\big(\Gamma_{s_i}\backslash\Omega\big)_{i=1}^n\nonumber \\
\;[\omega]_t & \longmapsto & \big([g_i^{-1}(\omega)]_{s_i}\big)_{i=1}^n,
\label{LastLine}
\end{eqnarray}
where $s_i\in S$ such that $[s_i] = [t_i]$ in 
$\GL(k)\backslash\GL(\Af)/\GL(\Ahat)$, and $g_i\in\GL(k)$ such that
$g_is_ia_i=t_i$ for some $a_i\in\GL(\Ahat)$, for $i=1,\ldots,n$. 

We note that $\GL(k)\backslash\GL(\Af)/\big(Z(\Af)\GL(\Ahat)\big) =
\GL(k)\backslash\GL(\Af)/\GL(\Ahat)$, as for any $x\in\Af$ we have 
$\big(\begin{smallmatrix}x & 0\\0 & x\end{smallmatrix}\big) \equiv 
\prod_{\fp}\left(\begin{smallmatrix}\pi_{\fp}^{n_{\fp}} & 0\\0 & 
\pi_{\fp}^{n_{\fp}}\end{smallmatrix}\right) \bmod\GL(\Ahat)$, where
$\pi_{\fp}\in k$ denotes a chosen uniformizer at $\fp$, and
$n_{\fp}=\min(0,v_{\fp}(x))$ is zero for almost all $\fp$. This latter scalar
is in $\GL(k)$, and, as scalars commute in $\GL(\Af)$, is killed by the left
$\GL(k)$ action. The same principle does not hold when we have $n>1$ copies of
$\GL(\Af)$ and of $Z(\Af)\GL(\Ahat)$ but only one copy of $\GL(k)$, which is why the
$Z(\Af)$ appears explicitly in the definition of~$\SC_n(\Cinf)^{\an}$.

Let 
\[
T^0 := \ker\big(T \rightarrow \GL(k)^n\backslash\GL(\Af)^n/\GL(\Ahat)^n\big).
\] 
Then for every $t\in T^0$, the Drinfeld modular curve $\Gamma_t\backslash\Omega$
is mapped into $\M^n(\Cinf)$ by $\theta$, where its image is a pure special
curve, as defined in the Introduction. 

Conversely, all pure special curves in $\M^n$ arise in this way:

\begin{proposition}\label{classification}
We have bijections
\[
\begin{array}{rcccp{5.5cm}}
T^0 & \longleftrightarrow & (Z(k)\GL(A))^n\backslash\GL(k)^n/\GL(k) &
\longleftrightarrow & $\{$Isomorphism classes of pure special curves in
$\M^n\}$\\
t & \longmapsto & [g] = [(g_1,\ldots,g_n)] & \longmapsto & Pure special curve 
defined by \mbox{$\omega\mapsto[(g_1(\omega),\ldots,g_n(\omega))]$,} 
\end{array}
\]
with $g\in\GL(k)^n$ such that $g=at^{-1}$ for some $a\in\big(Z(\Af)\GL(\Ahat)\big)^n$.
\end{proposition}

\proof In the first bijection, let $t\in T^0$, then by definition there exist
$g\in\GL(k)^n$ and $a^{-1}\in\big(Z(\Af)\GL(\Ahat)\big)^n$ such that $gta^{-1}=1$. Then
$[g]=[at^{-1}]$ is well-defined, since if also $g'=a't^{-1}$, then
$g'=a'a^{-1}g$ with $a'a^{-1}\in\big(Z(\Af)\GL(\Ahat)\big)^n\cap\GL(k)^n =
\big(Z(k)\GL(A)\big)^n$. Moreover, given $g\in\GL(k)^n$, there exist $h\in\GL(k)$,
$t\in T^0$ and  $a\in\big(Z(\Af)\GL(\Ahat)\big)^n$ such  that $g^{-1}=hta$, hence the
map is surjective. If also $g^{-1}=h't'a'$, then it is clear that $t=t'$, by
definition of $T$.

The second bijection follows as $(g_1,\ldots,g_n)$ and $(g'_1,\ldots,g'_n)$ in
$\GL(k)^n$ define the same pure special curve in $\M^n$ if and only if $g'_i =
\gamma_ig_i\sigma$ for all $i=1,\ldots,n$, for some
$(\gamma_1,\ldots,\gamma_n)\in\big(Z(k)\GL(A)\big)^n$ and $\sigma\in\GL(k)$. \qed

This proposition justifies our definition of $\SC_n$ as the space of all pure
special subcurves of $(\MM^2_{A,\Cinf})^n$. We also denote by
$\SC^0_n\subset\SC_n$  the subfamily of those Drinfeld modular curves
corresponding to $T^0$, i.e. $\SC^0_n(\Cinf)^{\an} \cong \coprod_{t\in
T^0}\Gamma_t\backslash\Omega$, which is the space of all pure special subcurves
of $\M^n$.

It remains to define {\em special subvarieties} of $\big(\MM^2_{A,\Cinf}\big)^n$
as products of pure special subcurves and CM points.

\subsection{Pure special curves and trees}

Let $\fp\subset A$ be a prime. Recall (e.g. \cite{SerreTrees}) that the {\em
Br\^uhat-Tits tree} $\TT_{\fp}$ of $\GL(k_{\fp})$ is the $(|\fp|+1)$-regular
tree whose vertices represent homothety classes of $A_{\fp}$-lattices in
$k_{\fp}^2$, and two vertices representing classes $\Lambda$ and $\Lambda'$ are
joined by an (unoriented) edge if there exist representative lattices 
$L\in\Lambda$ and $L'\in\Lambda'$ such that $L'\subset L$ and $L/L' \cong
A/\fp$. The group $\GL(k_{\fp})$ acts transitively on $\TT_{\fp}$ and the
stabilizer of the ``origin'' $v_{o,\fp}$, which is the vertex corresponding to
the lattice $A_{\fp}^2$, is $Z(k_{\fp})\GL(A_{\fp})$. Thus we get a bijection
$\TT_{\fp}\leftrightarrow\GL(k_{\fp})/Z(k_{\fp})\GL(A_{\fp})$.

Hence we have bijections 
\begin{equation}\label{tree1}
T \longleftrightarrow \GL(k)\backslash\GL(\Af)^n/\big(Z(\Af)\GL(\Ahat)\big)^n
\longleftrightarrow \GL(k)\backslash{\prod_{\fp}}'\TT_{\fp}^n.
\end{equation}
Here the restricted product $\prod'_{\fp}$ denotes the subset of those families
$(v_{1,\fp},\ldots,v_{n,\fp})_{\fp}$ of $n$-tuples of vertices
$v_{i,\fp}\in\TT_{\fp}$ such that $v_{1,\fp}=\cdots=v_{n,\fp}=v_{o,\fp}$ for
almost all primes $\fp$. 

We also have the following bijection
\begin{equation}\label{tree1.5}
\GL(k)\backslash {\prod_{\fp}}'\TT_{\fp}^n \longleftrightarrow
S\times {\prod_{\fp}}'\GL(k)\backslash\TT_{\fp}^n 
\quad \big( \leftrightarrow\Pic(A)\times 
{\prod_{\fp}}'\GL(k_{\fp})\backslash\TT_{\fp}^n\;\big)
\end{equation}
which arises from \pagebreak
\[
\GL(k)\backslash\GL(\Af)^n/\big(Z(\Af)\GL(\Ahat)\big)^n 
\stackrel{\sim}{\longrightarrow} 
\]
\[ \GL(k)\backslash\GL(\Af)/Z(\Af)\GL(\Ahat) \quad\times\quad 
\GL(\Af)\backslash\GL(\Af)^n/\big(Z(\Af)\GL(\Ahat)\big)^n; 
\]
\[
\;[g_1,\ldots,g_n] \longmapsto  
\big([g_1],[1,g_1^{-1}g_2,\ldots,g_1^{-1}g_n]\big),
\]
which is readily verified. This bijection may also be of help to the reader
puzzling over \cite[\S1.3]{BreuerPrep}.

The moduli interpretation of (\ref{tree1}) and (\ref{tree1.5}) is the following.

We may view $\TT_{\fp}$ as the tree of $\fp$-isogenies of Drinfeld modules as
in the elliptic curve case, but with the added subtlety that different vertices
will map to different irreducible components of $\big(\MM_{A,\Cinf}^2\big)^n$
when $\fp\subset A$ is not principal. Now let $t\in T$ correspond to the pure
special subcurve 
$C_t\stackrel{\theta}{\hookrightarrow}\big(\MM_{A,\Cinf}^2\big)^n$. Let
$(s,\underline{v})\in S\times {\prod_{\fp}}'\GL(k)\backslash\TT_{\fp}^n$
correspond to $t$ via (\ref{tree1.5}). We may choose a representative family
$(v_{1,\fp},\ldots,v_{n,\fp})_{\fp}$ of the class $\underline{v}$ such that
$v_{1,\fp}=v_{o,\fp}$ for all $\fp$. Then a typical point 
$x=(x_1,\ldots,x_n)\in C_t(\Cinf)$ satisfies: $x_1$ lies on the $s$-component
of $\MM_{A,\Cinf}^2$, and the isogenies $x_1\rightarrow
x_2\rightarrow\cdots\rightarrow x_n$ correspond to
$(v_{2,\fp},\ldots,v_{n,\fp})_{\fp}$, respectively.

In the case $n=3$ we have a particularly pleasing combinatorial description of
$\GL(k)\backslash\TT_{\fp}^3$. Any triple of vertices
$(v_1,v_2,v_3)\in\TT_{\fp}^3$ has a well-defined {\em center} $v_c$, defined by
the property that the three paths (possibly of length zero) from $v_c$ to the
$v_i$ are pairwise edge-disjoint. Denote by $(n_1,n_2,n_3)\in\NN_0^3$ the lengths of these
paths. Furthermore, $\GL(k_{\fp})$ acts $3$-transitively on the set of ends of
$\TT_{\fp}$, and it follows that a triple of vertices $(v_1,v_2,v_3)$ with
corresponding center $v_c$ and triple of path-lengths $(n_1,n_2,n_3)$ is mapped
by $\GL(k)$ to another triple $(v'_1,v'_2,v'_3)$ with center $v'_c$ and
path-lengths $(n'_1,n'_2,n'_3)$ if and only if
$(n_1,n_2,n_3)=(n'_1,n'_2,n'_3)$. Thus we have bijections 
\begin{equation}\label{tree2}
\GL(k)\backslash\TT_{\fp}^3 \longleftrightarrow \NN_0^3; \quad
{\prod_{\fp}}'\GL(k)\backslash\TT_{\fp}^3 \longleftrightarrow 
\II_A^3; \quad T \longleftrightarrow \Pic(A)\times\II_A^3,
\end{equation}
where $\II_A$ denotes the semigroup of non-zero $A$-ideals. The second bijection
is given by 
\[
(v_{1,\fp},v_{2,\fp},v_{3,\fp})_{\fp} \mapsto
(n_{1,\fp},n_{2,\fp},n_{3,\fp})_{\fp} \mapsto (\fn_1,\fn_2,\fn_3) = 
\prod_{\fp}\big(\fp^{n_{1,\fp}},\fp^{n_{2,\fp}},\fp^{n_{3,\fp}}\big).
\]

Let $Y$ be the special curve corresponding to the data 
$\big([\fa],(\fn_1,\fn_2,\fn_3)\big)\in \Pic(A)\times\II_A^3$. Then $Y$ maps
into $\M^3$ if and only if $[\fa]=1$ and $\fn_i=\langle N_i\rangle$ is
principal for $i=1,2,3$. In this case we have $p_{i,j}(Y) =
Y'_0(N_iN_j)\subset\M^2$. In particular, we have shown

\begin{lemma}\label{triples}
The set of pure special subcurves of $\M^3$ is in bijection with $(A/\Fq^\times)^3$.\qed
\end{lemma}

In general, if a pure modular curve
$Y\subset\big(\MM_{A,\Cinf}^2\big)^3$ corresponds to a triple
$(\fn_1,\fn_2,\fn_3)\in\II_A^3$, then $p_{i,j}(Y)$ is isomorphic to a suitable 
irreducible component of $\MM_A^2(\KK_0(\fn_i\fn_j))_{\Cinf}$.

The following result will be important later on.

\begin{proposition}\label{BigDegree}
Let $Y\subset\M^3$ be a pure special curve corresponding to a triple
$(\fn_1,\fn_2,\fn_3)$, and fix $(x_1,x_2)\in p_{1,2}(Y(\Cinf))$. Then
\[
|p_{1,2}^{-1}(x_1,x_2)\cap Y(\Cinf)| \geq 
\prod_{\fp|\fn_3}\frac{(|\fp|-1)(|\fp|+1)^{n_{3,\fp}-1}}{2n_{3,\fp}+1},
\]
where $\fn_i = \prod_{\fp}\fp^{n_{i,\fp}}$.
\end{proposition}

\proof Let $(v_{1,\fp},v_{2,\fp},v_{3,\fp})_{\fp}\in\prod'_{\fp}\TT_{\fp}^3$ be
a representative family of tuples corresponding to the special curve $Y$, and
denote by $(v_{c,\fp})_{\fp}$ its family of centers. Fixing $x_1$ and $x_2$
amounts to fixing $v_{1,\fp}$, $v_{2,\fp}$ for all $\fp$. We now count the
possible vertex families $(v_{3,\fp})_{\fp}$ for which
$(v_{1,\fp},v_{2,\fp},v_{3,\fp})_{\fp}$ represents $Y$. As
$(\fn_1,\fn_2,\fn_3)$ is prescribed, so is $v_{c,\fp}$ for every $\fp$, and
from (\ref{tree2}) follows that we must count the paths of length $n_{3,\fp}$
from $v_{c,\fp}$ which are edge-disjoint from the two paths leading to $v_{1,\fp}$
and $v_{2,\fp}$, which gives at least $(|\fp|-1)(|\fp|+1)^{n_{3,\fp}-1}$ for
each $\fp$. Thus the number of valid cyclic $\fn_3$-isogenies from $x_c$
(corresponding to $(v_{c,\fp})_{\fp}$) to $x_3$ is given by 
$\prod_{\fp|\fn_3}(|\fp|-1)(|\fp|+1)^{n_{3,\fp}-1}$. 

But if $x_c$ has complex multiplication by an order in $K$, then distinct
isogenies may well lead to the same $x_3$, they correspond to non-trivial
endomorphisms $f\in\End(x_c)$ of norm $N_{K/k}(f)=\fn_3^2$. There are at most
$\prod_{\fp}(2n_{3,\fp}+1)$ such endomorphisms, which completes the proof. 
\qed

\subsection{Main results and reduction to $\M^n$}

Let $\KK_i\subset\GL(\Ahat)$ be subgroups of finite index, for $i=1,\ldots,n$.
With a bit more effort, we could give a general treatment of pure special
subcurves of $\prod_{i=1}^n\MM_A^2(\KK_i)_{\Cinf}$ as in the previous two
sections. But the following definition is much easier:

Consider the morphism
\[
\pi : \prod_{i=1}^n\MM_A^2(\KK_i)_{\Cinf}\longrightarrow 
\big(\MM^2_{A,\Cinf}\big)^n
\]
induced by the inclusions $\KK_i\hookrightarrow\GL(\Ahat)$. Then an irreducible
subvariety $Y\subset \prod_{i=1}^n\MM_A^2(\KK_i)_{\Cinf}$ is a {\em special
subvariety} if and only if the image $\pi(Y)$ is a special subvariety of
$\big(\MM^2_{A,\Cinf}\big)^n$. Similarly for (pure) special subcurves.

We restate our main result as follows

\bigskip\noindent {\bf Theorem \ref{MainResult}$^\prime$\;\;}{\em
Let $\KK_i\subset\GL(\Ahat)$ be subgroups of finite index for $i=1,\ldots,n$.
Let $Y\subset\prod_{i=1}^n\MM_A^2(\KK_i)_{\Cinf}$ be an irreducible subvariety.
Then $Y(\Cinf)$ contains a Zariski-dense subset of CM points if and only if $Y$
is a special subvariety.}\bigskip

As every irreducible component of $\MM_A^2(\KK_i)_{\Cinf}$ is a Drinfeld
modular curve corresponding to $\Gamma_s\backslash\Omega$, for the arithmetic
subgroup $\Gamma_s = s\KK_i s^{-1}\cap \GL(k)$ with some $s\in\GL(\Af)$, we see
that the above Theorem is equivalent to Theorem \ref{MainResult}.

Next we want to show that it suffices to prove Theorem \ref{MainResult} for
$X=\M^n$.

\begin{proposition}\label{reduction}
Let $\Gamma'_i \subset \Gamma_i \subset \GL(k)$ be arithmetic subgroups
corresponding to Drinfeld modular curves $X'_i$ and $X_i$, respectively, for
$i=1,\ldots,n$. Let $X = X_1\times\cdots\times X_n$ and 
$X' = X'_1\times\cdots\times X'_n$. Then the canonical map 
$f : X' \rightarrow X$ preserves special subvarieties. In other words, a
subvariety $Y \subset X'$ is special if and only if $f(Y)$ is special.
\end{proposition}

\proof This is immediate, as special subvarieties are defined purely in
terms of isogeny relations between coordinates, and the property that any
constant projections $Y \rightarrow X'_i$ have CM points as images. These
properties are independent of any level structures. \qed

Now let $\Gamma_i\subset\GL(k)$ be an arithmetic subgroup such that
$X_i(\Cinf)^{\an} \cong \Gamma_i\backslash\Omega$, for each $i=1,\ldots,n$. Then
each $\Gamma'_i = \GL(A)\cap\Gamma_i$ is again an arithmetic subgroup, and we
have maps 
\[
\Gamma'_i\backslash\Omega \longrightarrow \Gamma_i\backslash\Omega, 
\quad \Gamma'_i\backslash\Omega \longrightarrow \GL(A)\backslash\Omega
\]
corresponding to morphisms of Drinfeld modular curves
\[
f_i : X'_i \longrightarrow  X_i, \quad g_i : X'_i \longrightarrow \M,
\]
for each $i=1,\ldots,n$. 

Applying Proposition \ref{reduction} to the maps $\prod_i f_i$ and $\prod_i g_i$
shows that Theorem \ref{MainResult} holds for $X_1\times\cdots\times X_n$ if 
and only if it holds for $\M^n$.

With this reduction step out of the way, the rest of our proof of Theorem
\ref{MainResult} will follow \cite{BreuerPrep} very closely, with $\M^n$
playing the role of $\AA^n$. The main geometric aspect of $\AA^n$ used in
\cite{BreuerPrep} is its structure as the product of affine lines, and it turns
out that the product structure of $\M^n$ will suffice, with some caveats.
Firstly, $\M$ is not defined over $k$ but over $H$, the Hilbert class field of
$k$. (Classically, the elliptic modular curve $Y(1)$ is defined over the Hilbert
class field of $\QQ$, which is just $\QQ$ itself. In our case, $k\neq H$ in
general, which illustrates the ``metamathematical'' phenomenon of splitting of
the two different roles of $\QQ$ into two different fields ($k$ and $H$) in
characteristic $p$. See \cite[Preface]{Goss} for a discussion).

Secondly, in general the isomorphism 
\[
\GL(A)\backslash\Omega \stackrel{\sim}{\longrightarrow} \M(\Cinf)^{\an}
\]
is no longer induced by a well-behaved ``$j$-invariant''  
$j : \Omega \rightarrow \Cinf$, as is the case for $A=\Fq[T]$. Consequently, we
cannot take advantage of the analytic properties of the $j$-function, and will
say nothing about the Weil heights of CM points.

\section{Hecke correspondences and CM points}
\reset

\subsection{The general formalism}

We briefly describe Hecke correspondences on $\MM_A^2(\KK)$ in general. Let
$g\in\GL(\Af)$, and set $\KK_g = \KK\cap g^{-1}\KK g$. Then $g$ acts from the
left on $\MM_A^2(\KK_g)$ (by letting $g^{-1}$ act from the {\em right} on
Drinfeld modules with full level structures, see \cite[II.3]{GekelerDMC}).
Combined with the standard projection $\pi : \MM_A^2(\KK_g) \rightarrow
\MM_A^2(\KK)$ this gives rise to the {\em Hecke correspondence} on
$\MM_A^2(\KK)$:
\begin{equation}\label{HeckeScheme}
T_g : \MM_A^2(\KK) \stackrel{\pi}{\longleftarrow} \MM_A^2(\KK_g)
\stackrel{\pi\circ g}{\longrightarrow} \MM_A^2(\KK).
\end{equation}
This is an algebraic correspondence of finite degree $[\KK : \KK_g]$. 

Analytically, this correspondence is described as follows. The group $\GL(\Af)$
acts from the left on $\GL(\Af)\times\Omega$ by $g\cdot (h,\omega) =
(hg^{-1},\omega)$. This gives rise to the correspondence
\[
\xymatrix{
\GL(\Af)\times\Omega\ar[r]^{g}\ar[d] & \GL(\Af)\times\Omega\ar[d]\\
\GL(k)\backslash\GL(\Af)\times\Omega/\KK & 
\GL(k)\backslash\GL(\Af)\times\Omega/\KK,
}
\]
which is easily seen to factor through $\GL(k)\times\KK_g$, thus describing the
finite correspondence
\begin{equation}
\xymatrix{
& \GL(k)\backslash\GL(\Af)\times\Omega/\KK_g\ar[dl]_{\pi}\ar[dr]^{\pi\circ g}
&\\
\GL(k)\backslash\GL(\Af)\times\Omega/\KK & T_g\ar@{.>}[l]\ar@{.>}[r] & 
\GL(k)\backslash\GL(\Af)\times\Omega/\KK.
}
\end{equation}

Explicitly, this correspondence acts as follows. Choose a finite set 
$\{g_i\in\GL(\Af) \;|\; i\in I_{g}\}$ of representatives for the cosets
$\KK\backslash\KK g\KK$. Let
$[h,\omega]\in\GL(k)\backslash\GL(\Af)\times\Omega/\KK \cong
\MM_A^2(\KK)(\Cinf)$. Then $T_g$ maps the point $[h,\omega]$ to the set 
$\{[hg_i^{-1},\omega] \;|\; i\in I_g\}.$

At first sight, the $g$-action does not seem to effect the $\Omega$ part of the
above spaces, but we will show next how $T_g$ acts on
$\MM_A^2(\KK)(\Cinf)\cong\coprod_{s\in S}\Gamma_s\backslash\Omega$. For
simplicity, we assume that $\det(\KK_g) = \det(\KK)$. That way,
$\MM_A^2(\KK_g)_{\Cinf}$ and $\MM_A^2(\KK)_{\Cinf}$ have the same number of
irreducible components, and we may choose the same set $1\in S\subset\GL(\Af)$
of representatives for $\GL(k)\backslash\GL(\Af)/\KK$ as for
$\GL(k)\backslash\GL(\Af)/\KK_g$. We may also assume that $h\in S$.

For each $i\in I_g$, there exist $f_{h,g,i}\in\GL(k)$ and $k_{h,g,i}\in\KK$
such that $f_{h,g,i}hg_i^{-1}k_{h,g,i}=s_{h,g} \in S$. Now $[hg_i^{-1},\omega]
= [s_{h,g}, f_{h,g,i}(\omega)]$, and this element maps to $[f_{h,g,i}(\omega)]
\in \Gamma_{s_{h,g}}\backslash\Omega$. Moreover, we can check that $s_{h,g}$
here only depends on $h$ and on $g$, not on $i$, so the correspondence $T_g$
maps the $h$-component $\Gamma_h\backslash\Omega$ to the $s_{h,g}$-component
$\Gamma_{s_{h,g}}\backslash\Omega$. Within the $h$-component, $T_g$ maps the point
$[\omega]$ to the set $\{[f_{h,g,i}(\omega)]\;|\; i\in I_g\}$, and in fact the
set $\{f_{h,g,i} \;|\; i\in I_g\}$ is a set of representatives for
$\Gamma_{s_{h,g}}\backslash\Gamma_{s_{h,g}}f_{h,g}\Gamma_h$, where
$f_{h,g}\in\GL(k)$ is such that $f_{h,g}hg^{-1}k_{h,g}=s_{h,g}\in S$ for some
$k_{h,g}\in\KK$. Hence, for each $h\in S$, the correspondence $T_g$ on
$\MM_A^2(\KK)$ induces the following correspondence from the $h$-component to
the $s_{h,g}$-component of $\MM_A^2(\KK)(\CC)$:
\[
\xymatrix{
\Omega\ar[d]\ar[r]^{f_{h,g}} & \Omega\ar[d] \\
\Gamma_h\backslash\Omega & \Gamma_{s_{h,g}}\backslash\Omega
}
\]
which factors through 
\[
\Gamma_h\cap f_{h,g}^{-1}\Gamma_{s_{h,g}}f_{h,g} = h\KK_gh^{-1}\cap \GL(k)
= \Gamma_{h}\cap\Gamma_{hg},
\]
thus giving the correspondence
\begin{equation}\label{HeckeAnalytic}
\xymatrix{
& (\Gamma_h\cap\Gamma_{hg})\backslash\Omega\ar[dl]_{\pi}\ar[dr]^{\pi\circ
f_{h,g}} & \\
\Gamma_h\backslash\Omega & & \Gamma_{s_{h,g}}\backslash\Omega
}
\end{equation}
which has finite degree $[\Gamma_h : \Gamma_h\cap\Gamma_{hg}] = [\KK : \KK_g]$.
Thus the correspondence (\ref{HeckeAnalytic}) gives a component by component
description of the correspondence $T_g$ on $\MM_A^2(\KK)_{\Cinf}$ of 
(\ref{HeckeScheme}).

\Remark Let $s_1,s_2 \in S$, and set $g=s_2^{-1}s_1$. Then
$f_{s_1,g}=k_{s_1,g}=1$, and so $T_g$ maps the $s_1$-component of
$\MM_A^2(\KK)_{\Cinf}$ to the $s_2$-component, while leaving the $\omega$
unchanged in that component. 

We also notice that $s_{h,g}=h$ if and only if $g\in\KK h^{-1}\GL(k)h$.

\subsection{Hecke correspondences on $\M^n$}

In the remainder of this section, we closely follow \cite[\S2]{BreuerPrep}.
From now on, we will assume that $\KK=\GL(\Ahat)$ and $g=\diag(N,1,\ldots,1)$
for some $N\in A$ with $|N|>1$. In this case, we see that $\KK_g =
\KK_0(\langle N \rangle) =: \KK_0(N)$, $\Gamma_1=\GL(A)$,
$\Gamma_1\cap\Gamma_g=\Gamma_0(N)$, and furthermore, $T_N:=T_g$ acts trivially
on the set of irreducible components on $\MM^2_{A,\Cinf}$, in particular, it
induces a correspondence, again denoted $T_N$, on $\M$. The correspondence is
also described by the inclusion $Y'_0(N) \subset\M^2$ from (\ref{Y'}), and it
maps the point $x\in\M(\Cinf)$ to the set of those points $y\in\M(\Cinf)$
corresponding to Drinfeld modules linked to $x$ via cyclic $\langle
N\rangle$-isogenies.

We define $T_{\M^n,N}$ to be the correspondence on $\M^n$ which is the product
of the correspondences $T_N$ on each factor $\M$. When there is no risk of
confusion, we also denote it by $T_N$. We say $T_{\M^n,N}$ {\em stabilizes} an
algebraic subvariety $Y\subset\M^n$ if $Y\subset T_{\M^n,N}(Y)$. In this case,
we define the {\em restriction of $T_N$ to $Y$}, denoted $T_{Y,N}$, to be the
union of the irreducible components of $T_{\M^n,N}\cap (Y\times Y)$ of maximal
dimension. It is a correspondence on $Y$ which is still surjective, in the
sense that the two projections $T_{Y,N} \rightarrow Y$ are surjective. Whenever
we mention $T_{Y,N}$, it is implicit that $T_N$ stabilizes $Y$. 

Let $x\in\M(\Cinf)$, and suppose that $x\in T_N(x)$. Then $x$ admits a cyclic
endomorphism, hence is a CM point. For a given $N\in A,\; |N|>1$, there are only
finitely many points stabilized by $T_N$, which correspond to the points of the
diagonal in $\M^2$ which intersect $Y'_0(N)$. On the other hand, given a CM
point $x\in\M(\Cinf)$ with $R=\End(x)$, there are infinitely many $N\in A$ such
that $x\in T_N(x)$, namely all those $N$ composed of primes $\fp\subset A$ for
which $\fp R$ is a product of two distinct principal prime ideals of $R$.
Equivalently, these primes split completely in the class field corresponding to
$\Pic(R)$, hence the set of such primes has density $1/2|\Pic(R)|$, by
\v{C}ebotarev. Notice also that such primes $\fp$ are necessarily principal in
$A$.

\subsection{Some intersection theory}

$\M$ is affine, so we may fix an embedding $\M\subset\AA^m_{\Cinf}$. Then we
obtain embeddings $\M^n\subset\AA^{mn}_{\Cinf}$ and
$\overline{\M^n}\subset\PP^{mn}_{\Cinf}$. (Here $\bar{\cdot}$ denotes the
Zariski-closure).

For any irreducible $Y\subset\M^n$, we denote by $\deg(Y)$ the degree of 
$\overline{Y}\subset\PP^{mn}$ in the usual sense. If $Y=\cup_i Y_i$ is a union
of irreducible components, then $\deg(Y) := \sum_i\deg(Y_i)$, (regardless of
their dimensions).

We denote by $\psi(N) = [\GL(A):\Gamma_0(N)] = \prod_{\fp|N}(1+|\fp|^{-1})$ the
degree of the correspondence $T_N$ on $\M$. In particular, 
$\psi(N)\leq\deg(Y'_0(N))\leq 2\psi(N)$.

We collect the following facts:

\begin{proposition}\label{intersections}
Let $Y\subset\M^n$ be an algebraic subvariety.
\begin{enumerate}
\item $Y$ has at most $\deg(Y)$ irreducible components.
\item If $Y'\subset\M^n$ is another subvariety, then $\deg(Y\cap Y')
\leq\deg(Y)\deg(Y')$.
\item $\deg(T_{\M^n,N}(Y)) \leq 2^n\psi(N)^n\deg(Y)$.
\item There are only finitely many pure special subcurves of $\M^n$ of degree
less than a given bound.
\end{enumerate}
\end{proposition}

\proof (1) is trivial, (2) and (3) follow from \cite[8.4.6]{Fulton}, and  (4)
follows from the fact that a curve $Y\subset\M^n$ is pure special if and only
if $p_{i,j}(Y) = Y'_0(N_{ij})$ for some $N_{ij}\in A$ for all $i\neq j$, and
$\deg(Y'_0(N_{ij}))\geq\psi(N_{ij})$. \qed

\subsection{Preimages in $\Omega^n$}

Denote by $\pi$ the rigid analytic map $\pi :
\Omega^n\rightarrow\M^n(\Cinf)^{\an}$, which is the quotient for the
$\GL(A)^n$-action. For each irreducible component $Y_i$ of $Y$, we choose an
irreducible rigid analytic variety $Z_i\subset\Omega^n$ with
$\pi(Z_i)=Y_i(\Cinf)^{\an}$. The group $\GL(\kinf)^n$ acts on $\Omega^n$, and
the $\GL(A)^n$-orbit of $Z_i$ is $\pi^{-1}(Y_i)$. We next describe the action
of $T_N$ in $\Omega^n$.

We may choose a set $T\subset\GL(k)$ of representatives for
$\GL(A)\backslash\GL(A)\big(\begin{smallmatrix}N & 0 \\ 0 & 1
\end{smallmatrix}\big)\GL(A)$ such that the following holds: For every pair of
principal ideals $\fa,\fd\subset A$ with $\fa\fd=\langle N \rangle$, and set of
representatives $\{b_1,\ldots,b_r\} \subset A$ of $A/\fd$, there exist $a,d\in
A$ with $\fa=\langle a \rangle$ and $\fd=\langle d \rangle$ such that
\begin{equation}\label{Tstr}
\begin{pmatrix}a & b_j \\ 0 & d\end{pmatrix} \in T\quad\text{for
$j=1,\ldots,r$.} 
\end{equation}
One easily verifies that these elements indeed represent distinct coset classes.
We write $T^n=\{t_j\in\GL(k)^n \;|\; j\in\II^n\}$, where 
$\II=\{1,\ldots,\psi(N)=|T|\}$. 

For each $Z_i$ we define $\JJ_{Z_i}\subset\II^n$ as the set of those indices
$j$ for which $t_j(Z_i)\subset\pi^{-1}(Y_i)$. Now let $y\in Y_i(\Cinf)$ and 
choose some $z\in Z_i$ with $\pi(z)=y$. Then
\begin{eqnarray}
T_{\M^n,N}(y) & = & \{\pi(t_j(z)) \;|\; j\in\II^n\}, \\
T_{Y,N}(y) & = & \{\pi(t_j(z)) \;|\; j\in\JJ_{Z_i}\}.
\end{eqnarray}
In particular, we see that each $\JJ_{Z_i}$ is non-empty, as $T_{Y,N}$ is a
surjective correspondence. In fact, under suitable conditions, the sets
$\JJ_{Z_i}$ are fairly large:

\begin{proposition}\label{surjectivity}
Let $Y\subset\M^n$ be a subvariety all of whose irreducible components have the
same dimension, and suppose that $Y\subset T_{\M^n,N}(Y)$ for some $N\in A$ such
that $\langle N \rangle$ is a product of distinct primes $\fp\subset A$ of even
degree satisfying $|\fp|\geq\max(13,\deg(Y))$. 
Let $I\subset\{1,\ldots,n\}$ and let $Y_i$ be an irreducible
component of $Y$ for which the projection $p_I : Y_i\rightarrow\M^I$ is
dominant, and choose a preimage
$Z_i\subset\Omega^n$ of $Y_i(\Cinf)^{\an}$. Then the projection
\[
p_I : \JJ_{Z_i} \longrightarrow \II^I
\]
is surjective.
\end{proposition}

\proof This is \cite[Thm. 4]{BreuerPrep}, the proof is exactly the same. We
briefly explain the hypotheses on $N$: we use the fact that $\M$ has a
$\PSL(A/NA)\cong\prod_{\fp}\PSL(A/\fp)$-covering $Y_2(N)$, of which $Y_0(N)$ is
a subcover (Proposition \ref{covering}), and that this group has no proper
subgroups of index $\deg(Y)$ or less when $|\fp|\geq\max(13,\deg(Y))$.\qed

\subsection{An interlude in group theory}

We remark that the $\GL(\kinf)$ action on $\Omega$ induces a $\PGL(\kinf)$
action, as the center acts trivially. Until now we have found it more convenient
to work with $\GL$, but in order to continue we will need a number of
group-theoretic results, and here it will be simpler to work with $\PGL$, as had
been done throughout \cite{BreuerPrep}.

For the convenience of the reader, we recall here some basic properties of the
groups $\PGLki$, which we will need. These results were used implicitly in
\cite{BreuerPrep}, and thus may also help the reader with that paper.

\begin{lemma}
$\kinf^{\times}/\kinf^{\times 2} \cong (\ZZ/2\ZZ)^2$.
\end{lemma}

\proof This is trivial, as $\kinf\cong\Fq((\varpi))$ for a uniformizer
$\varpi$. \qed

\begin{proposition}\label{Group1}
Every non-trivial normal subgroup of $\PGLki$ contains $\PSLki$. In particular,
$\PGLki$ has no non-trivial discrete normal subgroups.
\end{proposition}

\proof Let $H \triangleleft \PGLki$ be a normal subgroup. Then, as $\PSLki$ is
simple, either $\PSLki \subset H$, or $H\cap\PSLki = \{1\}$. In the latter case
we get an embedding  $H \hookrightarrow \PGLki/\PSLki \cong
\kinf^{\times}/\kinf^{\times 2}  \cong (\ZZ/2\ZZ)^2$. It remains to show that
$\PGLki$ has no normal subgroups isomorphic to $\ZZ/2\ZZ$ or $(\ZZ/2\ZZ)^2$.
This may be verified with explicit calculations, by conjugating elements of
order $2$ with  $\big(\begin{smallmatrix} 1 & 1 \\ 0 &
1\end{smallmatrix}\big)$. \qed

\begin{proposition}\label{Group2}
Let $H$ be a subgroup of finite index in $G=\PGLki$. Then $H$ is normal and 
contains $\PSLki$. In particular, if $H$ is simple then $H=\PSLki$.
\end{proposition}

\proof $G$ acts on the cosets $G/H$, giving a representation  $\rho : G
\rightarrow \Aut(G/H) \cong S_n$ where $n = [G:H]$. Then  $K=\ker(\rho)$ is a
normal subgroup of $G$ contained in $H$. As $K$ is infinite, $K$ must intersect
non-trivially with $\PSLki$, hence contains $\PSLki$, as the latter is simple.
We will not need the fact that $H$ is normal, but it is known that all the
subgroups between $\PSLki$ and $\PGLki$ are normal in $G$. \qed

\begin{corollary}\label{Group3}
Let $\PSLki \subset H \subset \PGLki$ and suppose that 
$f : H \hookrightarrow \PGLki$ is a monomorphism with image of finite index.
Then $f|_{\PSLki}$ in an automorphism of $\PSLki$.
\end{corollary}

\proof The image $f(\PSLki)$ has finite index in $\PGLki$ and is simple, hence
$f(\PSLki)=\PSLki$. \qed

The above results may easily be generalized to $\PGL$ over other (infinite) 
fields. Mostly we have just used the fact that 
$\kinf^{\times}/\kinf^{\times 2}$ is finite.

\begin{proposition}[Goursat's Lemma]\label{Goursat}
Let $G_1$ and $G_2$ be groups, and $H\subset G_1\times G_2$ a subgroup such that
the two projections $pr_i : H\rightarrow G_i$ are surjective. Then
$K_i=\ker(pr_i)$ can be considered a normal subgroup of $G_j$, for $i\neq j$,
and $H$ is the inverse image of the graph of an isomorphism $\rho : G_1/K_2
\stackrel{\sim}\rightarrow G_2/K_1$.
\end{proposition}

\proof This is straight forward. The map
\begin{eqnarray*}
\rho:G_1/K_2 & \longrightarrow & G_2/K_1 \\
         g_1 & \longmapsto & \text{$g_2$ with $(g_1,g_2)\in H$}
\end{eqnarray*}
is easily checked to be a well-defined isomorphism. Now $H$ is the inverse
image of the graph of $\rho$.
\qed

\subsection{Curves stabilized by Hecke correspondences}

Our next goal is to characterize special subvarieties of $\M^n$ by their
property of being stabilized by suitable Hecke correspondences. The hard part
is to prove this for curves.

\begin{theorem}\label{TopCurves} 
Let $Y\subset\M^2$ be an irreducible algebraic curve such that both projections
$Y\rightarrow\M$ are dominant, and suppose $Y\subset T_{\M^2,N}(Y)$ for some
$N\in A$ such that $\langle N \rangle$ is a product of distinct primes
$\fp\subset A$ of even degree satisfying $|\fp|\geq\max(13,\deg(Y))$. Then
$Y=Y'_0(N')$ for some $N'\in A$. 
\end{theorem}

Note that we cannot deduce $N'$ from $N$. Indeed $Y'_0(N')$ is stabilized by
$T_{\M^2,N}$ for all $N\in A$ coprime to $N'$.\bigskip

\proof Suppose that the hypotheses of Theorem \ref{TopCurves} are satisfied.
The group $G:=\PGLki^2$ acts on $\Omega^2$, and we define the following
subgroups: $S:=\PSLki^2$, $\Gamma:=\PGL(A)^2$ and $\Sigma:=\PSL(A)^2$. We fix
an irreducible rigid analytic curve $Z\subset\Omega^2$ with
$\pi(Z)=Y(\Cinf)^{\an}$. Let $G_Z\subset G$ be the stabilizer of $Z$, which is
a closed analytic subgroup of $G$. We also define the subgroups $S_Z = G_Z\cap
S$, $\Gamma_Z = G_Z\cap\Gamma$ and $\Sigma_Z  = G_Z\cap\Sigma$. We intend to
prove Theorem \ref{TopCurves} by investigating the structure of $G_Z$ and
$S_Z$.

Denote by $pr_i : G_Z \rightarrow \PGLki$ the two projections for $i=1,2$. 

\begin{lemma}\label{GroupProj}
The projection $pr_i : G_Z \rightarrow \PGLki$ is injective, and 
$pr_i(\Gamma_Z)$ has finite index in $\PGL(A)$, for $i=1,2$.
\end{lemma}

\proof Exactly the same as \cite[Lemma 2.11]{BreuerPrep}. Here one uses 
Proposition \ref{Group1}. \qed

Let $i\in\JJ_Z$, then $t_i(Z)\subset\pi^{-1}(Y)$, so there is some
$\gamma_i\in\Gamma$ such that $\gamma_i t_i\in G_Z$. As the projections
$pr_1,pr_2 : \JJ_Z\rightarrow\II$ are surjective (Proposition
\ref{surjectivity}), we get many non-trivial elements of $G_Z$ this way. 

Let $H_1 = pr_1(G_Z)$. Our first goal is to show that $H_1$ contains $\PSLki$. By a
slight abuse of notation, we also denote the elements of $pr_1(T)\subset\GL(k)$
by $t_i$, for $i\in\II$. We have seen that, for every $i\in pr_1(\JJ_Z)=\II$,
there exists some $\gamma_i\in\PGL(A)$ such that $g_i:=\gamma_i t_i \in H_1$.
Lemma \ref{GroupProj} says that $H_1\cap\PGL(A)$ has finite index in $\PGL(A)$,
and we choose a finite set $R\subset\PGL(A)$ of representatives for
$\PGL(A)/(\PGL(A)\cap H_1)$. 

We claim that given any string $i_1,\ldots,i_n$ of elements of $\II$, and any
$a\in\GL(A)$, there exists $\gamma\in R$ such that 
\[
\gamma t_{i_n}t_{i_{n-1}}\cdots t_{i_1}a \in H_1.
\]
Indeed, by induction it suffices to prove the claim for $n=1$. Let $i_1\in\II$
and $a\in\GL(A)$ be given. $\GL(A)$ acts from the right on the set of left
cosets $\GL(A)\backslash\GL(A)\big(\begin{smallmatrix}N & 0\\ 0 &
1\end{smallmatrix}\big)\GL(A)$, so we let $a$ act from the right on
$\GL(A)\cdot t_{i_1}$, obtaining $\GL(A)\cdot t_{i_1}a = \GL(A)\cdot t_j$ for
some $j\in\II$. Thus $t_{i_1}a = \gamma'_jt_j$, and for suitable $\gamma\in R$ 
and $\gamma'\in H_1\cap\PGL(A)$ we have
\[
\gamma t_{i_1}a = \gamma\gamma'_jt_j = \gamma(\gamma'_j\gamma_j^{-1})\gamma_jt_j
= \gamma'\gamma_jt_j = \gamma'g_j\in H_1.
\]
This proves the claim.

Now, multiplying by a suitable power of $\big(\begin{smallmatrix}N & 0 \\ 0 &
N\end{smallmatrix}\big)$, we see from (\ref{Tstr}) that for any $x\in A[1/N]$
and $a\in\GL(A)$, there exists some $\gamma_{x,a}\in R$ such that 
$\gamma_{x,a}\big(\begin{smallmatrix}1 & x\\ 0 & 1\end{smallmatrix}\big)a \in
H_1$.  Denote by $E\subset\PSL(A[1/N])$ the subgroup generated by elementary
matrices. We have shown that $H_1\cap E$ has finite index in $E$. As $A[1/N]$
is dense in $\kinf$, and $\PSLki$ is generated by elementary matrices, it
follows that $E$ is dense in $\PSLki$. 

Next, we see that $H_1$ is a closed subgroup of $\PGLki$, exactly the same way
as \cite[Lemma 2.12 and the following paragraph]{BreuerPrep}, where we need
only replace $\PSL(A[1/\fm])$ by $E$. It follows that $H_1\cap\PSLki$ has
finite index in $\PSLki$, so by Proposition \ref{Group2} 
\begin{equation}\label{containsPSL}
\PSLki\subset H_1,\quad\text{and similarly, $\PSLki\subset H_2=pr_2(G_Z)$.} 
\end{equation}

Now, since the projections $pr_1,pr_2 : G_Z\rightarrow\PGLki$ are injective, 
Proposition \ref{Goursat} implies that 
\[
G_Z = \{(g,\rho(g)) \;|\; g\in H_1\},
\]
where $\rho : H_1 \stackrel{\sim}{\rightarrow} H_2$ is an isomorphism. From
(\ref{containsPSL}) and Corollary \ref{Group3} follows that
\begin{equation}\label{S_Z}
S_Z = G_Z\cap\PSLki^2 = \{(g,\rho(g)) \;|\; g\in\PSLki\},
\end{equation}
where $\rho|_{\PSLki}$ is an automorphism of $\PSLki$.

It is known that the automorphisms of $\PSLki$ are all of the form $g\mapsto
hg^{\sigma}h^{-1}$ for some $h\in\PGLki$ and $\sigma\in\Aut(\kinf)$, see
\cite{Hua}. By the definition of $\Sigma_Z$ and (\ref{S_Z}), we see that
$h\cdot pr_1(\Sigma_Z)^{\sigma}\cdot h^{-1}\subset\PSL(A)$. On the other hand,
Lemma \ref{GroupProj} tells us that $pr_1(\Sigma_Z)$ has finite index in
$\PSL(A)$. This in turn tells us a lot about $h$ and $\sigma$. Fix some
$\varpi\in A$ so that $\varpi^{-1}$ is a uniformizer for $\kinf$, i.e.
$\kinf=\Fq((\varpi^{-1}))$.

\begin{proposition}\label{Autom}
Let $G$ be a subgroup of finite index in $\PGL(A)$, and suppose that 
$hG^{\sigma}h^{-1}\subset\PGL(k)$, for some $h\in\PGLki$ and
$\sigma\in\Aut(\kinf)$. Then $h\in\PGL(k)$, $\sigma|_{\Fq}\in\Gal(\Fq/\FF_p)$,
and $\sigma(\varpi^{-1}) =  u\varpi^{-1}+v$ for some $u\in\Fq^*$ and $v\in\Fq$.
\end{proposition}

\proof Let $R=\Fq[[\varpi^{-1}]]$ be the unique valuation ring in $\kinf$.
We show that $h\in\PGL(k)$, $\sigma(k)\subset k$ and $\sigma(R)\subset R$
exactly as in \cite[Prop. 2.13]{BreuerPrep}. 

It follows that $\sigma(A)\subset A$ and $\sigma(A^{\times})\subset
A^{\times}$, hence $\sigma(\Fq)\subset\Fq$. Furthermore,
$\sigma(\varpi^{-1})\in A^{-1}$ must again be a uniformizer for $\kinf$. But
all other uniformizers of $\kinf$ that lie in $A^{-1}$ are of the form
$u\varpi^{-1}+v$ for some $u\in\Fq^*$ and $v\in\Fq$. This completes the proof.
\qed 

Now that we have assembled enough ingredients, the proof of Theorem 
\ref{TopCurves} follows exactly as in \cite[\S2.7]{BreuerPrep}. We provide a
sketch. We have seen that
\[
S_Z = \{(g,hg^{\sigma}h^{-1}) \;|\; g\in\PSLki\},
\]
where $h\in\PGL(k)$ and $\sigma\in\Aut(\kinf)$ is as in Proposition
\ref{Autom}. In particular, there exists an integer $t\geq 0$ such that
$\sigma(\alpha)=\alpha^{q^t}$, for all $\alpha\in\Fq$.  
Let $f=(\varpi^{-q}-\varpi^{-1})^{q-1}$,
and $F:=\FF_p((f))$, which is a complete subfield of $\kinf$ on which $\sigma$
acts trivially. Fix some non-square $\alpha\in\Fq$, and set
\[
P = \{z\in\Omega \;|\; z^2=\alpha e,\; e\in F\}\subset\Omega.
\]
Pick any $z_1\in P$. Then $S_1=\Stab_{\PSL(F)}(z_1)$ is a one-dimensional Lie
group (so far we have made essential use of the assumption that $p$ is odd).

Now let $z_2\in\Omega$ such that $(z_1,z_2)\in Z$. Then the ``$S_1$-orbit'' 
\[
\{(g(z_1),hg^{\sigma}h^{-1}(z_2)) \;|\; g\in S_1\}\subset
Z\cap(\{z_1\}\times\Omega) 
\]
is discrete, but $S_1$ is not, so there exists $1\neq g\in S_1$ such that
$hg^{\sigma}h^{-1}$ fixes both $z_2$ and $h'(z_1)$, where $h'=h\circ
\left(\begin{smallmatrix}\alpha^{(p^t-1)/2} & 0\\ 0 & 1\end{smallmatrix}\right)
\in\PGL(k)$.
This means that $z_2$ and $h'(z_1)$ are conjugate over $\kinf$, i.e. $z_2 =\pm
h'(z_1)$. It follows that either $\pi(z_1,h'(z_1))\in Y(\Cinf)$ or
$\pi(-z_1,h'(-z_1))\in Y(\Cinf)$. Both of these points also lie on $Y'_0(N')$,
where $N'=\deg(h')$ is independent of $(z_1,z_2)$. As $P$ is uncountable,
whereas the fibers of $\pi$ are countable, it follows that $Y=Y'_0(N')$.
\qed

To extend Theorem \ref{TopCurves} to subvarieties $Y\subset\M^n$ of higher
dimension we may follow \cite[\S2.8 and Corollaries 2.8 and 2.9]{BreuerPrep}
almost verbatim, no new ingredients are required. We must just keep in mind
that $\M$ is defined over $H$, the Hilbert class field of $k$. We also use
Galois action on CM points as in the following section.

\begin{theorem}\label{TopVar}
Let $F$ be a field lying between $H$ and $\Cinf$. Let $Y\subset\M^n$ be an
$F$-irreducible algebraic subvariety, containing a CM point $x\in Y(\Cinf)$.
Suppose that $Y\subset T_{\M^n,N}(Y)$ for some $N\in A$ such that $\langle N
\rangle$ is a product of distinct primes $\fp\subset A$ of even degree
satisfying $|\fp|\geq\max(13,\deg(Y))$. Then $Y\subset\M^n$ is a special
subvariety. \qed
\end{theorem}

\subsection{CM points on curves}

We may now start proving our main results, by exploiting the behavior of CM
points under Galois action (Proposition \ref{CMmain}) in conjunction with
Theorems \ref{TopCurves} and \ref{TopVar}. We first treat the case of curves,
where our results are effective.

\begin{theorem}\label{MainCurves}
Let $X=X_1\times\cdots\times X_n$ be a product of Drinfeld modular curves. Let
$F/H$ be a finite extension, and $d\in\NN$. Then there exists an absolutely
computable constant $B=B(X,F,d)>0$ such that the following holds. Let $Y\subset
X$ be an irreducible algebraic subcurve of degree $d$ and defined over $F$. Then
$Y$ is a special subcurve if and only if $Y(\Cinf)$ contains a CM point $x$
satisfying $\HCM(x)>B$.
\end{theorem}

This indeed implies Theorem \ref{MainResult} for curves, since if $Y\subset X$
contains a Zariski-dense set of CM points, then it contains CM points of
arbitrary CM height and $Y$ is defined over some finite extension $F/k$, as the
CM points are all defined over $k^\sep$.\bigskip

\proof We again follow \cite[\S3.4]{BreuerPrep} very closely, but will provide
full details here for the benefit of the reader.

It follows from Proposition \ref{reduction} that we may assume that 
$X=\M^n$. Furthermore, we may assume that none of the projections $p_i :
Y\rightarrow\M$ are constant. Using the fact that $Y\subset\M^n$ is a pure
special subcurve if and only if $p_{i,j}(Y)\subset\M^2$ is special for all
$i<j$, we have reduced the problem to the case $n=2$.

Let $x=(x_1,x_2)\in Y(\Cinf)$ be a CM point. For $i=1,2$ we write 
$R_i=\End(x_i)=A+\ff_i\OO_{K_i}$, an order of conductor $\ff_i\subset A$ in the
CM field $K_i$, which has genus $g_i$. Furthermore, we denote by
$K_i(x_i):=H_{R_i}$ the ring class field of $R_i$, which is a field of
definition for $x_i$, and we have $\Gal(K_i(x_i)/K_i)\cong\Pic(R_i)$, and
$K_i(x_i)/K_i$ is unramified outside $\ff_i$. We denote by $K=K_1K_2$ and
$K(x_1,x_2)=K_1(x_1)K_2(x_2)$ the composite fields. Denote by $F_s$ the
separable closure of $k$ in $F$ (which contains $H$), and by $L$ the Galois
closure of $F_sK(x_1,x_2)$ over $k$.

Let $\fp$ be a prime of $k$ of even degree which splits completely in $F_sK$
(in particular, $\fp$ is principal, as it splits in $H$), and suppose
$\fp\nmid\ff_1\ff_2$. Let $\fP$ be a prime of $L$ lying above $\fp$, and denote
by $\fP_i$ its restriction to $K_i(x_i)$. Denote by $\sigma\in\Aut(FL/FK)$ an
extension of the Frobenius element $(\fP,L/k)$, and let
$\sigma_i=\sigma|_{K_i(x_i)}=(\fP_i,K_i(x_i)/K_i)$ (remember that $\fp$ splits
in $K$ and is unramified in $L$).

It follows from Proposition \ref{CMmain} that $x_i$ and $x_i^{\sigma}$
correspond to Drinfeld modules linked by cyclic $\fp$-isogenies, and hence
\begin{equation}
(x_1,x_2)\in Y\cap T_{\M^2,\fp}(Y^{\sigma}) = Y\cap T_{\M^2,\fp}(Y).
\end{equation}

Moreover, the whole $\Gal(FK(x_1,x_2)/F)$-orbit of $(x_1,x_2)$ lies in this
intersection, which thus contains at least
$\max\big(|\Pic(R_1)|,|\Pic(R_2)|\big)/[F:k]$ points. On the other hand, from
Proposition \ref{intersections}, $\deg(Y\cap T_{\M^2,\fp}(Y)) \leq
4(|\fp|+1)^2\deg(Y)^2$. Therefore, if
\begin{equation}\label{bigPic}
|\Pic(R_i)|/[F:k] > 4(|\fp|+1)^2\deg(Y)^2, \quad\text{for $i=1$ or $2$,} 
\end{equation}
then the intersection is improper, $Y\subset T_{\M^2,\fp}(Y)$ and hence $Y$ is
special (Theorem \ref{TopCurves}), provided also $|\fp|\geq\max(13,\deg(Y))$.

It remains to show that such a suitable prime $\fp$ indeed exists, if $\HCM(x)$
is sufficiently large. Let $M$ be the Galois closure of $F_sK$ over $k$, and set
\[
\pi_M(t) = \#\{\fp\subset A \;|\; \text{prime, split in $M$ and $|\fp|=q^t$}\}.
\]
Let $T\in k$ be a transcendental element such that $k$ is a finite separable
geometric extension of $\Fq(T)$, and let $e=[k:\Fq(T)]$. Let $\FF$ be the
algebraic closure of $\Fq$ in $M$, let $n_c=[\FF:\Fq]$ be the constant
extension degree and $n_g=[M:\FF k]$ the geometric extension degree of $M/k$.
The \v{C}ebotarev Theorem for function fields \cite[Prop. 5.16]{FriedJarden}
says
\begin{equation}
\text{if $n_c|t$, then}\quad |\pi_M(t)-\frac{1}{n_g}q^t/t| < 
  4\big(e^2+g_M(e+1)/2 +g_k+1\big)q^{t/2},
\end{equation}
where $g_M$ and $g_k$ are the genuses of $M$ and $k$, respectively. We may bound
$g_M$ in terms of $g_1,g_2$ and $g_k$ using the Castelnuovo inequality
\cite[III.10.3]{Stichtenoth}, and eventually obtain
\begin{equation}
\pi_M(t) > C_1q^t/t - \big(C_2(g_1+g_2)+C_3\big)q^{t/2},
\end{equation}
where $C_1,C_2$ and $C_3$ are absolutely computable positive constants,
depending on $k$ and $F$.

We want $\pi_M(t)>\log_q|\ff_1\ff_2|$, $q^t \geq\max(13,\deg(Y))$ and
$2n_c|t$, so that there exists a prime $\fp$ which splits in $M$ (and thus in
$F_sK$), does not divide $\ff_1\ff_2$, and satisfies the hypotheses of Theorem
\ref{TopCurves}. We also want (\ref{bigPic}) to hold, for which we employ
Proposition \ref{ClassNumber}: 
\[
|\Pic(R_i)|>C_{\varepsilon}\HCM(x_i)^{1-\varepsilon} =
C_{\varepsilon}(q^{g_i}|\ff_i|)^{1-\varepsilon},\quad\text{for any
$\varepsilon>0$.}
\]

In summary, we need a simultaneous solution $t\in 2n_c\NN$ to the three 
inequalities
\[
q^t \geq\max(13,\deg(Y)),
\]
\[
C_1q^t/t - \big(C_2(g_1+g_2)+C_3\big)q^{t/2} > \log_q|\ff_1\ff_2|,
\]
and
\[
C_{\varepsilon}(q^{g_i}|\ff_i|)^{1-\varepsilon} > 4[F:k](q^t+1)^2\deg(Y)^2
\quad\text{for some $\varepsilon>0$, and $i=1$ or $2$.}
\]
Such a solution will always exist if $\HCM(x)=\max(\HCM(x_1),\HCM(x_2))$ is
sufficiently large. (Intuitively, $q^t$ must be large compared to
$\log_q|\ff_i|$ and $g_i$, and small compared to $|\ff_i|$ and $q^{g_i}$).
\qed

\subsection{Completing the proof of Theorem \ref{MainResult}}

The proof of Theorem \ref{MainResult} now follows exactly as in
\cite[\S3.5]{BreuerPrep}. We sketch the proof here for the sake of completeness.

Firstly, as $Y(\Cinf)$ contains a Zariski-dense subset $S$ of CM points, which
are defined over $k^{\sep}$, there exists a finite Galois extension $F/k$ over
which $Y$ is defined. We will use induction on $d=\dim(Y)$, the case $d=1$
following from Theorem \ref{MainCurves}, so we assume that $d\geq 2$ and $n\geq
3$. We may furthermore assume that $Y\subset\M^n$ is a hypersurface, as it is an
irreducible component of 
\[
\bigcap_{\substack{I\subset\{1,\ldots,n\}\\ |I|=d+1}}p_I^{-1}p_I(Y).
\]
Lastly, we may assume that all the projections $p_i : Y\rightarrow\M$ are
dominant.

\paragraph{Step 1.} For a given constant $B>0$, we may assume that every
$x=(x_1,\ldots,x_n)\in S$ satisfies $\HCM(x_i)>B,\;\forall i=1,\ldots,n$
(otherwise we replace $S$ by a Zariski-dense subset). 

Pick one such $x\in S$. {\bf Suppose} that there exist primes
$\fp_1,\ldots,\fp_{d-1}\subset A$ of even degree satisfying the following
conditions:
\begin{enumerate}
\item[\em (i)] Each $\fp_j$ splits in $F$ and in $\End(x_i)$, for all $i=1,\ldots,n$.
\item[\em (ii)] $|\fp_1|\geq\max(13,\deg(Y))$.
\item[\em (iii)] $|\fp_{j+1}|\geq (\deg(Y))^{2^j}\prod_{m=1}^j(2|\fp_m|+2)^{n2^{j-m}}$, for
$j=1,\ldots,d-2$.
\item[\em (iv)] $|\Pic(\End(x_i))| > [F:k]|\fp_{d-1}|^2(2|\fp_{d-1}|+2)^n$, for all
$i=1,\ldots,n$.
\end{enumerate}
Then one shows, again using Galois action on $x$ together with Theorem
\ref{TopVar}, that there exists a pure special subvariety $Y_x\subset Y$ containing $x$.

As the points in $S$ are Zariski-dense, it follows that there exists a
Zariski-dense family $\C$ of pure special subcurves $C\subset Y$, $C\in\C$.

\paragraph{Step 2.} We now show that $Y$ is special. Choose a CM point
$x_1\in\M(\Cinf)$, and consider the slice
\[
Y_1 = Y\cap(\{x_1\}\times\M^{n-1}).
\]
Each pure special curve $C\in\C$ intersects $Y_1$ in at least one CM point, and
we denote by $Y'$ the Zariski-closure of these intersection points:
\[
Y' = \overline{\cup_{C\in\C}(C\cap Y_1)} \subset Y_1.
\]
We must have $\dim(Y')<\dim(Y)$, so by the induction hypothesis, $Y'$ is
special. We write $Y'=Y'_1\cup\cdots\cup Y'_r$ as a union of irreducible
components. Replacing $\C$ by a Zariski-dense subfamily and renumbering if
necessary, we may assume that $Y'_1$ contains at least $1/r$ of the points of
$C\cap Y_1$ for all $C\in\C$.

Now, either $Y'_1\cong \{y\}\times\M^m$ for some CM point
$y\in\M^{n-m}(\Cinf)$, in which case $Y=Y'_0(N')\times\M^{n-2}$ is special
(here we use the fact that Theorem \ref{MainResult} has already been proved for curves), or
else at least one pure special curve appears as a factor of $Y'_1$. 

It follows that there exist indices $1<i<j$ such that $p_{i,j}(Y'_1)=Y'_0(M)$
for some {\em fixed} $M\in A$. Fix $C\in\C$ and let the pure special curve
$p_{\{1,i,j\}}(C)\subset\M^3$ correspond to the triple
$(N_{C,1},N_{C,i},N_{C,j})\in A^3$ via Lemma \ref{triples}. Again restricting $\C$ and switching $i$ and $j$ if
necessary, we may assume that $|N_{c,i}|\leq |N_{C,j}|$ for all $C\in\C$.

Now we fix $C\in\C$ and $x_1,x_i\in\M(\Cinf)$. Then the number of distinct
$x_j\in\M(\Cinf)$ such that $(x_1,x_i,x_j)\in p_{\{1,i,j\}}(C)$ is bounded from
below by an increasing function in $|N_{C,j}|$ (Proposition \ref{BigDegree}).
But at least $1/r$ of these points $x_j$ must also satisfy $(x_i,x_j)\in
Y'_0(M)$, of which there can be at most $\psi(M)$. It follows that $|N_{C,i}|$
and $|N_{C,j}|$ are bounded independently of $C\in\C$. Restricting $\C$ once
again, we may assume that $p_{i,j}(C)=Y'_0(N_0)$ for all $C\in\C$. 

Now one can show that $Y\cong Y'_0(N_0)\times\M^{n-2}$, which is special. 

\paragraph{Step 3.} It remains to show that the primes
$\fp_1,\ldots,\fp_{d-1}\subset A$ satisfying {\em (i)}-{\em (iv)} above actually exist, if
the constant $B>0$ is chosen sufficiently large. Let $x=(x_1,\ldots,x_n)\in S$
such that $\HCM(x_i)=q^{g_i}|\ff_i|>B$ for all $i=1,\ldots,n$. 

As before, the problem boils down to finding simultaneous solutions
$t_1,\ldots,t_{d-1}\in 2n_c\NN$ to the following four inequalities:
\[
q^{t_1}\geq\max(13,\deg(Y)),
\]
\[
C_{\varepsilon}(q^{g_i}|\ff_i|)^{1-\varepsilon} >
[F:k]q^{2t_{d-1}}(2q^{t_{d-1}}+2)^n,\quad\text{for some $\varepsilon>0$ and some
$1\leq i\leq n$,}
\]
\[
q^{t_{j+1}} \geq (\deg(Y))^{2^j}\prod_{m=1}^j(2q^{t_m}+2)^{n2^{j-m}},\quad
\text{for all $1\leq j \leq d-1$,}
\]
\[
C_1q^{t_j}/t_j - \big(C_2(g_1+\cdots+g_n)+C_3\big)q^{t_j/2} >
\log_{q}|\ff_1\cdots\ff_n|,\quad\text{for all $1\leq j \leq d-1$.}
\]
If $B>0$ is sufficiently large, then such solutions exist. This completes the
proof of Theorem~\ref{MainResult}. \qed


\section{Application to Heegner points}
\label{HeegnerSection}
\reset

In this section we apply Theorem \ref{MainResult} to extend the main result of
\cite{Breuer04} to arbitrary global function fields of odd characteristic.

Let $E$ be an elliptic curve defined over $k$, which we recall is a global
function field over $\Fq$, with $q$ odd. 
Suppose that the $j$-invariant of $E$ is not constant, so we say that $E$ is non-isotrivial. Then, replacing
$k$ by a finite extension if necessary, there exists a place of $k$ at which 
$E$ has split multiplicative reduction, and if we call this place $\infty$ we 
are in the situation of the previous sections.
Now the conductor of $E$ is of the form 
$\fn\cdot\infty$, for an ideal $\fn\subset A$. It follows from the work of 
Drinfeld and others that we have a modular parametrization
\begin{equation}\label{ModParam}
X_0(\fn) \longrightarrow E
\end{equation}
defined over $k$, where $X_0(\fn)$ is the smooth projective model for the curve
$Y_0(\fn)$ defined in \S\ref{Y0section}. See \cite{GekelerReversat} for a detailed treatment.


We fix a prime $\fp$ of A for the remainder of this section.

\begin{lemma}
\label{Heegner}
There exist infinitely many quadratic imaginary extensions $K/k$ satisfying the following two conditions:
\begin{itemize}
\item[(i)] Every prime $\fq\subset A$, $\fq\not=\fp$, which ramifies in $K/k$ is principal in $k$.
\item[(ii)] Every prime $\fq\subset A$ which divides $\fn$ splits in $K/k$ {\em (Heegner hypothesis)}. 
\end{itemize}
\end{lemma}

\proof
Denote by $k_\fn$ the ray class field of $k$ with conductor $\fn$. Then a prime $\fq\subset A$ splits completely in $k_\fn$ if and only if $\fq = \langle x \rangle$ with $x\equiv 1 \bmod\fn$. 
Denote by $\cQ_\fn$ the set of primes $\fq\subset A$ of odd degree which split completely in $k_\fn$.
By the \v{C}ebotarev Theorem \cite[Prop. 5.16]{FriedJarden}, this set is infinite. 
%
%
Now let $m\in A$ such that $\langle m \rangle$ is a product of primes in $\cQ_\fn$ and $\deg(m)$ is odd. Then $k(\sqrt{m})/k$ is a quadratic imaginary extension satisfying conditions {\em (i)} and {\em (ii)} above.
\qed

Fix a quadratic imaginary extension $K/k$ satisfying conditions \ref{Heegner}.{\em (i)}-{\em (ii)} above, 
and let $n\in\NN$.
Denote by $\OO_K$ the ring of integers of $K$, and let $\OO_n = A + \fp^n\OO_K$,
which is an order of conductor $\fp^n$ in $\OO_K$. Thanks to condition \ref{Heegner}.{\em (ii)},
there exists an ideal $\fN_n\subset\OO_n$ such that 
$\OO_n/\fN_n \cong A/\fn$ as $A$-modules. It follows that the pair of lattices $(\OO_n, \fN^{-1}_n)$ defines a pair of Drinfeld modules with complex multiplication by $\OO_n$ and linked by a cyclic $\fn$-isogeny. Thus the pair defines a point $x_n \in X_0(\fn)(K[\fp^n])$, where $K[\fp^n]$ denotes the ring class field of $\OO_n$.

We let $K[\infty]:=\cup_{n\geq 0}K[\fp^n]$ in a chosen algebraic closure $\bar{k}$ of $k$. Then 
$\Gal(K[\infty]/K)\cong G_0\times\ZZ_p^{\infty}$, where $\ZZ_p^{\infty}$ denotes the direct product of countably many copies of $\ZZ_p^+$, where $p$ is the characteristic of $k$, and $G_0$ is a finite abelian group, see \cite[Proposition 2.1]{Breuer04}. We denote by $H[\infty]\subset K[\infty]$ the fixed field of $G_0$, so $\Gal(K[\infty]/H[\infty])=G_0$ and $\Gal(H[\infty]/K)\cong\ZZ_p^{\infty}$. We write $H[\fp^n] = H[\infty]\cap K[\fp^n]$.

We now define the $n$th {\em higher Heegner point} on $E$ by
\begin{equation}\label{HeegnerPoint}
y_n := \Tr_{G_0}(\pi(x_n)) = \sum_{\sigma\in G_0}\pi(x_n^{\sigma})\in E(H[\fp^n]).
\end{equation}

The main result of this section is 
\begin{theorem}\label{HeegnerTheorem}
Let $I\subset\NN$ be an infinite subset. In the above situation, the group generated by 
$\{y_n \;|\; n\in I\}$ in $E(H[\infty])$ has finite torsion and infinite rank.
\end{theorem}

\proof
Most of the work has already been done in \cite[Theorem 2]{Breuer04}, combined with Theorem \ref{MainResult}. It remains to verify the surjectivity of a new modular parametrization $\pi' : X_0(\fm\fn) \rightarrow E$, which we proceed to describe.

Let $\fp_1,\ldots,\fp_g$ denote the primes $\not=\fp$ of $A$ which ramify in $K/k$. As we have assumed condition \ref{Heegner}.{\em (i)}, we know that they are all principal ideals. Let $\fm=\fp_1\cdots\fp_g$ denote their product. 
Choose a set $S\subset\GL(\Af)$ of simultaneous representatives for the double quotients 
$\GL(k)\backslash\GL(\Af)/\KK_0(\fn)$ and $\GL(k)\backslash\GL(\Af)/\KK_0(\fn\fm)$, (see \S\ref{Y0section} for notation). For each $s\in S$ we write $\Gamma_s(\fn)=s\KK_0(\fn)s^{-1}\cap\GL(k)$ and 
$\Gamma_s(\fn\fm)=s\KK_0(\fn\fm)s^{-1}\cap\GL(k)$. Setting $\Omega^*=\Omega\cup\PP^1(k)$, on which $\GL(k)$ acts in the obvious way, we get
\begin{eqnarray*}
X_0(\fn)(\Cinf)^{\an} & \cong & \coprod_{s\in S}\Gamma_s(\fn)\backslash\Omega^*\\
X_0(\fn\fm)(\Cinf)^{\an} & \cong & \coprod_{s\in S}\Gamma_s(\fn\fm)\backslash\Omega^*.
\end{eqnarray*}
Every divisor $\fd|\fm$ is principal, and we write $\fd=\langle d\rangle$ for a chosen $d\in A$. Let $D$ denote a set of these $d$'s as $\fd$ ranges through all divisors of $\fm$, so $|D|=2^g$. 

We define the {\em full degeneracy map} $\delta : X_0(\fn\fm) \rightarrow X_0(\fn)^{2^g}$ by its action on $\Cinf$-valued points: 
\begin{equation}\label{degen}
\delta : [\omega] \longmapsto \big([d\omega]\big)_{d\in D}\quad\mbox{on each $\Gamma_0(\fn\fm)\backslash\Omega^*$}.
\end{equation}

Next, we give an analytic description of the modular parametrization (\ref{ModParam}), see \cite{GekelerReversat} for details. Denote by $\TT_{\infty}$ the Br\^uhat-Tits tree of $\GL(\kinf)$, and by ${\underline{H\!}\,}(\TT_{\infty},\ZZ)$ the group of $\ZZ$-valued harmonic cochains on the set of edges of $\TT_{\infty}$. For each $s\in S$ one associates to $E$ a primitive Hecke newform $\varphi_s \in {\underline{H\!}\,}_{!}(\TT_{\infty},\ZZ)^{\Gamma_s(\fn)}$, the latter group denoting those harmonic cochains invariant under $\Gamma_s(\fn)$-action and with compact (= finite) support on $\Gamma_s(\fn)\backslash\TT_{\infty}$. To $\varphi_s$ one associates a certain holomorphic theta function 
$u_s : \Omega \rightarrow \Cinf^\times$ with multiplier $c_s : \Gamma_s(\fn)\rightarrow \Cinf^\times$; in other words, $u_s(\alpha\omega)=c_s(\alpha)u_s(\omega)$ for all $\omega\in\Omega$ and $\alpha\in\Gamma_s(\fn)$. We let 
$\Delta_s = \{c_s(\alpha) \;|\; \alpha\in\Gamma_s(\fn)\}$, which is a multiplicative lattice in $\Cinf^\times$. The elliptic curve $E$, which has split multiplicative reduction at $\infty$, is isomorphic to the Tate curve $\Cinf^\times/\Delta_s$, for each $s\in S$. Finally, on each $s$-component of $Y_0(\fn)$, the modular parametrization (\ref{ModParam}) is given explicitly on $\Cinf$-valued points by
\begin{eqnarray*}
\Gamma_s\backslash\Omega & \longrightarrow & \Cinf^\times/\Delta_s\\
\,[\omega] & \longmapsto & u_s(\omega)\;\bmod\Delta_s,
\end{eqnarray*}
and the cusps $X_0(\fn)\smallsetminus Y_0(\fn)$ map to the identity of $E$.

Now we can combine (\ref{ModParam}) with (\ref{degen}) to obtain a new modular parametrization of $E$:
\begin{eqnarray}
\pi' : X_0(\fn\fm) & \arrover{\Sigma\circ\pi\circ\delta} & E \\
\,[\omega] & \longmapsto & \prod_{d\in D}u_s(d\omega)\;\bmod\Delta_s, \quad\mbox{on each $\Gamma_s(\fn\fm)\backslash\Omega$}.\nonumber
\end{eqnarray}

According to \cite[Theorem 2]{Breuer04} it remains for us to show that 
$\pi' : \Gamma_s(\fn\fm)\backslash\Omega^* \rightarrow E$ is surjective for every\footnote{In \cite{Breuer04} we mistakenly only required $\pi' : X_0(\fn\fm)\rightarrow E$ to be surjective, whereas we actually need surjectivity for each irreducible component of $X_0(\fn\fm)$.} $s\in S$, in other words, that 
$u'_s(\omega) := \prod_{d\in D}u_s(d\omega)$ is not constant. 

Denote by $\OO_{\Omega}(\Omega)$ the ring of rigid analytic functions on $\Omega$, then there is an exact sequence
\[
1 \longrightarrow \Cinf^\times  \longrightarrow \OO_{\Omega}(\Omega)^\times \stackrel{r}{\longrightarrow}
{\underline{H\!}\,}(\TT_{\infty},\ZZ) \longrightarrow 0, 
\]
and the homomorphism $r$ satisfies $r(f\circ\alpha) = r(f)\circ\alpha$ for all $\alpha\in\GL(k)$. Moreover, 
$r(u_s) = \varphi_s$.

Now suppose that $u'_s$ is constant. Then 
\begin{equation}\label{ZeroSum}
0 = r(u'_s) = \sum_{d\in D}\varphi_s\circ\big(\begin{smallmatrix}d & 0 \\ 0 & 1\end{smallmatrix}\big).
\end{equation}
Now $\varphi_s\not=0$, so there exists an edge $e_0$ of $\TT_{\infty}$ such that $\varphi_s(e_0)\not=0$. 
Thus by (\ref{ZeroSum}) there is some $1\not=d_0\in D$ such that 
$\varphi_s\left(\big(\begin{smallmatrix}d_0&0\\0&1\end{smallmatrix}\big)\cdot e_0\right)\not=0$. 
Write $e_1=\big(\begin{smallmatrix}d_0&0\\0&1\end{smallmatrix}\big)\cdot e_0$,
and we again find some $1\not=d_1\in D$ such that 
$\varphi_s\left(\big(\begin{smallmatrix}d_1&0\\0&1\end{smallmatrix}\big)\cdot e_1\right)\not=0$, 
and so on, giving us an infinite sequence $e_0,e_1,\ldots,$ of edges of $\TT_{\infty}$ with 
$\varphi_s(e_i)\not=0$ for all $i\geq 0$. Moreover, as $|d_i|>1$ for all $i$, we see that the $e_i$'s 
form a sequence of distinct edges all lying on one end of $\TT_{\infty}$. It follows that their images in
$\Gamma_s(\fn)\backslash\TT_{\infty}$ still form an infinite sequence of edges lying on one end of 
$\Gamma_s(\fn)\backslash\TT_{\infty}$. But this contradicts
the fact that $\varphi_s$ has compact support modulo $\Gamma_s(\fn)$, which completes the proof of Theorem \ref{HeegnerTheorem}. \qed

\Remark The assumption \ref{Heegner}.{\em (i)} should not be necessary, but
the trouble begins when, if we have non-principal ramified primes other than $\fp$, the generalization of
(\ref{ZeroSum}) mixes Hecke newforms $\varphi_s$ for several different $s\in S$ in a single equation. 

\paragraph{Acknowledgment.} The author would like to thank Gebhard B\"ockle for 
pointing out the trick used to complete the proof of Theorem \ref{HeegnerTheorem}.


{\small
\begin{tabular}{p{9.5cm}}
\hline\\
Department of Mathematics, University of Stellenbosch, Stellenbosch 7600, South
Africa.\\
fbreuer@sun.ac.za
\end{tabular}}


\begin{thebibliography}{99}

\bibitem{Andre} Y.~Andr\'e, ``Finitude des couples d'invariants modulaires
    singuliers sur une courbe alg\'ebrique plane non modulaire'', {\em J. reine
    angew. Math.} {\bf 505} (1998), 203-208.

\bibitem{Breuer04}F.~Breuer, ``Higher Heegner points on elliptic curves over
    function fields'', {\em J. Number Theory} {\bf 104} (2004), 315-326.

\bibitem{BreuerPrep}F.~Breuer, ``The Andr\'e-Oort conjecture for products of
    Drinfeld modular curves'', to appear in {\em J. reine angew. Math.}
    ({\tt http://arxiv.org/abs/math.NT/0303038})
    
\bibitem{Edixhoven}S.J.~Edixhoven, ``Special points on products of modular
    curves'', to appear in {\em Duke Math. J.}
    ({\tt http://arxiv.org/abs/math.NT/0302138})

\bibitem{FriedJarden}M.~Fried and M.~Jarden, ``Field Arithmetic'',
    Springer-Verlag, 1986.

\bibitem{Fulton}W.~Fulton, ``Intersection Theory'', Springer-Verlag, 1984.

\bibitem{GekelerDMC} E.-U.~Gekeler, ``Drinfeld Modular Curves'', {\em Lecture
    Notes in Mathematics} {\bf 1231}, Springer-Verlag, 1986.
    
\bibitem{GekelerReversat} E.-U.~Gekeler and M.~Reversat, ``Jacobians of 
    Drinfeld modular curves'', {\em J. reine angew. Math.} {\bf 476} 
    (1996), 27-93.    

\bibitem{Goss}D.~Goss, ``Basic Structures of Function Field Arithmetic'',
    Springer-Verlag, 1996.

\bibitem{Hayes79} D.~Hayes, ``Explicit class field theory in global function
    fields'', in: ``Studies in algebra and number theory'' (G.C.Rota, ed.),
    Academic Press, New York, 1979.
    
\bibitem{Hayes92} D.~Hayes, ``A Brief introduction to Drinfeld modules'', in:
    The Arithmetic of Function Fields (eds. D.~Goss et al), de Gruyter, New
    York-Berlin, 1992.
    
\bibitem{Hua} L.K.~Hua, appendix to: J.~Dieudonn\'e, ``On the automorphisms of
    the classical groups'', {\em Memoirs Amer. Math. Soc.} {\bf 2} (1951), 1-95.

\bibitem{vdHeidenThesis}G.-J.~van der Heiden, ``Weil pairing and the Drinfeld
    modular curve'', Ph.D. thesis, Rijksuniversiteit Groningen, 2003.

\bibitem{Neukirch}J.~Neukirch, ``Algebraische Zahlentheorie'', Springer-Verlag,
    1992.

\bibitem{SerreTrees}J.-P.~Serre, ``Trees'', Springer-Verlag, 1980.

\bibitem{Stichtenoth}H.~Stichtenoth, ``Algebraic Function Fields and Codes'',
    Springer-Verlag, 1993.

\end{thebibliography}
\end{document}